\documentclass[10pt,a4paper,reqno]{amsart}%
\usepackage{amssymb}
\usepackage{amsmath}
\usepackage{amsfonts}
\usepackage[T1]{fontenc}
\usepackage{color}
\usepackage[usenames,dvipsnames]{xcolor}
\usepackage{bbm}
\usepackage{graphicx}
\usepackage{rotating}
\usepackage{hyperref}
\usepackage{mathrsfs}
\usepackage[all]{xy}%
\providecommand{\U}[1]{\protect\rule{.1in}{.1in}}
\newtheorem{theorem}{Theorem}

\newtheorem{example}[theorem]{Example}

\newtheorem{remark}[theorem]{Remark}

\thanks{}
\email{lvitagliano@unisa.it}
\begin{document}
\title[Characteristics and Singularities of Solutions of PDEs]{Characteristics, Bicharacteristics, and Geometric Singularities of Solutions
of PDEs}
\author{Luca Vitagliano}
\address{DipMat, Universit\`a degli Studi di Salerno, {\& Istituto Nazionale di Fisica
Nucleare, GC Salerno,} Via Ponte don Melillo, 84084 Fisciano (SA), Italy.}

\begin{abstract}
Many physical systems are described by partial differential equations (PDEs).
Determinism then requires the Cauchy problem to be well-posed. Even when the
Cauchy problem is well-posed for generic Cauchy data, there may exist
characteristic Cauchy data. Characteristics of PDEs play an important role
both in Mathematics and in Physics. I will review the theory of
characteristics and bicharacteristics of PDEs, with a special emphasis on
intrinsic aspects, i.e., those aspects which are invariant under general
changes of coordinates. After a basically analytic introduction, I will pass
to a modern, geometric point of view, presenting characteristics within the
jet space approach to PDEs. In particular, I will discuss the relationship
between characteristics and singularities of solutions and observe that:
\textquotedblleft wave-fronts are characteristic surfaces and propagate along
bicharacteristics\textquotedblright. This remark may be understood as a
mathematical formulation of the wave/particle duality in optics and/or quantum
mechanics. The content of the paper reflects the three hour minicourse that I
gave at the XXII International Fall Workshop on Geometry and Physics,
September 2--5, 2013, \'{E}vora, Portugal.

\end{abstract}

\maketitle

\tableofcontents

\section*{Introduction}

Many physical systems are described by partial differential equations (PDEs).
Determinism then requires the Cauchy problem to be well-posed. Even when the
Cauchy problem is well-posed for generic Cauchy surfaces, there may exist
\emph{characteristic Cauchy data}. Roughly speaking, characteristic Cauchy
data are those for which the Cauchy problem is ill-posed, in the sense of
non-existence or non-uniqueness of corresponding solutions. Surprisingly
enough, characteristic Cauchy data play an important role both in the
(mathematical) theory of PDEs and in Theoretical Physics. From a mathematical
point of view, characteristics of PDEs are related to intermediate integrals,
classifications of PDEs, singularities of solutions (besides Cauchy problems).
From a physical point of view, if one interprets independent variables as
space-time coordinates and dependent variables as fields, then a
characteristic Cauchy surface may be understood as the wave-front of a
\textquotedblleft bounded\textquotedblright\ disturbance in the fields,
propagating in the space-time. Characteristic Cauchy surfaces are often
themselves described by a first order, scalar PDE. In their turn first order,
scalar PDEs can be integrated with the method of characteristics. Namely, the
integration problem can be reduced to the integration problem for a system of
ordinary differential equations (ODEs) whose solutions \textquotedblleft
foliate\textquotedblright\ solutions of the original PDE. Accordingly,
characteristic surfaces are foliated by lines: \emph{characteristic lines} in
Cauchy terminology, \emph{bicharacteristic lines} in Hadamard terminology.
From a physical point of view, one concludes that a \emph{wave-front
propagates along bicharacteristics}. Notice that the transition between the
three different \textquotedblleft mathematical regimes\textquotedblright\
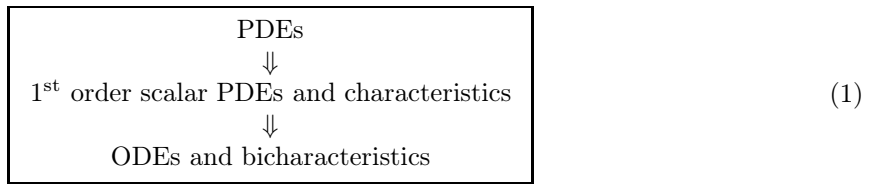
\begin{equation}
\fbox{$%
\begin{array}
[c]{c}%
\text{PDEs}\\
\Downarrow\\
\text{1}^{\text{st}}\text{ order scalar PDEs and characteristics}\\
\Downarrow\\
\text{ODEs and bicharacteristics}%
\end{array}
$} \label{13}%
\end{equation}
formalizes in rigorous terms the transition between three different
\textquotedblleft physical regimes\textquotedblright:
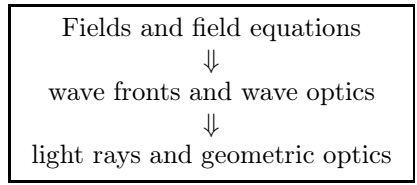
\[
\fbox{$%
\begin{array}
[c]{c}%
\text{Fields and field equations}\\
\Downarrow\\
\text{wave fronts and wave optics}\\
\Downarrow\\
\text{light rays and geometric optics}%
\end{array}
$}%
\]
Even more, the equation for characteristic surfaces is often an
Hamilton-Jacobi equation. It is well known that the Hamilton-Jacobi equation
is the short wave-length limit of the Schr\"{o}dinger equation. Actually,
interpreting the Hamilton-Jacobi equation as an equation for the wave-front of
a wave-function propagating in the space-time (see, e.g., \cite{o01}), one can
infer the Schr\"{o}dinger equation according to the analogy:
\[
\fbox{wave optics/geometric optics\quad=\quad wave mechanics/classical
mechanics}%
\]
In this sense transition (\ref{13}) is also analogous to the transition from
quantum mechanics to classical mechanics, summarized in the scheme:
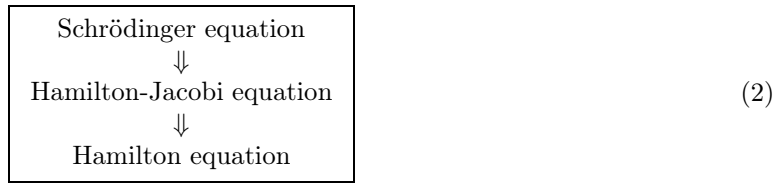
\begin{equation}
\fbox{$%
\begin{array}
[c]{c}%
\text{Schr\"{o}dinger equation}\\
\Downarrow\\
\text{Hamilton-Jacobi equation}\\
\Downarrow\\
\text{Hamilton equation}%
\end{array}
$} \label{14}%
\end{equation}
Accordingly, the \emph{quantizing}, i.e., reversing the arrows in
(\ref{14}),\emph{ }is analogous to \textquotedblleft\emph{reconstructing a PDE
from its (bi)characteristics}\textquotedblright. For certain specific classes
of PDEs the reconstruction can be actually accomplished, and, in a sense,
\emph{quantization is not ambiguous}. The aim of this paper is reviewing the
theory of (bi)characteristics of PDEs and its physical interpretation. In
particular, I will describe in some details the transition (\ref{13}) focusing
on intrinsic aspects, i.e., those aspects which are independent of the choice
of coordinates. Differential geometry will be then the natural language.

The paper is divided into three sections. In the first section, I discuss
Cauchy problems and characteristic Cauchy data. I conclude with some examples
from Mathematical Physics. This section is basically analytic and makes use of
local coordinates. However, most of the results therein are actually
independent of the choice of coordinates. In the second section, I present the
geometric setting for PDEs and their characteristics, specifically, jet
spaces. Characteristics of PDEs has a nice, intrinsic definition in terms of
jets. The geometric setting clarifies the relationship between characteristics
and singularities of solutions. In the last section, I focus on
bicharacteristics. Often characteristic surfaces are governed by a first order
scalar PDE $\mathcal{E}$. The geometry underlying such PDEs is contact
geometry which is at the basis of the method of characteristics. It may happen
that $\mathcal{E}$ is an Hamilton-Jacobi equation. There is a symplectic
version of the method of characteristics for Hamilton-Jacobi equations based
on the Hamilton-Jacobi theorem. This motives me to review the Hamilton-Jacobi
theory. I conclude speculating about the possibility of extending the
Hamilton-Jacobi theory to field theory in a covariant way, thus opening the
road through a \emph{rigorous, covariant, Schr\"{o}dinger quantization of
gauge theories}.

\section{Characteristic Cauchy Data for PDEs\label{Sec1}}

\subsection{Cauchy Problems}

The evolution of many physical systems, especially (but not only) in classical
physics, is described by a system of (sometimes non-linear) partial
differential equations (PDEs). Determinism requires that the full evolution of
the system is anambiguosly determined by the initial configuration. From a
mathematical point of view this means that the Cauchy problem for the
corresponding PDE should be well-posed, i.e., there should be
(\emph{existence}) a unique (\emph{uniqueness}) solution for any set of
(physically admissible) Cauchy data. The most general way to understand a set
of Cauchy data is \textquotedblleft a general hypersurface $\Sigma$ in the
space of independent variables + derivatives of the dependent variables normal
to $\Sigma$ along $\Sigma$ itself\textquotedblright. Even if the Cauchy
problem is well-posed for generic Cauchy data, existence or uniqueness may
fail for special Cauchy data usually referred as \emph{characteristic Cauchy
data}. Nonetheless, characteristic Cauchy data have a nice physical
interpretation. In this section I will recall some basic facts about the
Cauchy problem, characteristics of PDEs and their physical interpretation. For
simplicity, I will mainly focus on the case of determined systems of
quasi-linear partial differential equations. I will use local coordinates
everywhere, and I will conclude with few examples, mainly from Mathematical Physics.

\subsubsection{Cauchy problems in normal form}

Let $\boldsymbol{u}=(u^{1},\ldots,u^{m})$ be a vector valued function of the
$n$ real variables $x=(x^{1},\ldots,x^{n})$. I will often ineterpret the $x$'s
as space-time coordinates, and the $\boldsymbol{u}$'s as components of a field
propagating on the space-time. Put
\[
\boldsymbol{u}_{I}:=\dfrac{\partial^{|I|}}{\partial x^{I}}\boldsymbol{u}%
=\frac{\partial^{\ell}}{\partial x^{i_{1}}\cdots\partial x^{i_{\ell}}%
}\boldsymbol{u},
\]
where $I=i_{1}\cdots i_{\ell}$ is a multi-index (denoting multiple partial
derivatives) and $|I|{}:=\ell$ (the number of derivatives). Sometimes, it is
convenient to split space-time coordinates into \textquotedblleft space
coordinates\textquotedblright\ $\boldsymbol{x}=(x^{1},\ldots,x^{n-1})$ + a
time coordinate $t=x^{n}$. In this case I use the following notation for
multiple space-time derivatives:
\[
\boldsymbol{u}_{\ell,J}:=\dfrac{\partial^{|J|+\ell}}{\partial\boldsymbol{x}%
^{J}\partial t^{\ell}}\boldsymbol{u}.
\]
Consider the system of $m$ PDEs in $m$ unknown functions
\begin{equation}
\dfrac{\partial^{k}\boldsymbol{u}}{\partial t^{k}}=\boldsymbol{f}%
(t,\boldsymbol{x},\ldots,\boldsymbol{u}_{\ell,J},\ldots)\quad\label{N}%
\end{equation}
and the initial data problem
\begin{equation}
\left\{
\begin{array}
[c]{l}%
\dfrac{\partial^{k}\boldsymbol{u}}{\partial t^{k}}=\boldsymbol{f}%
(t,\boldsymbol{x},\ldots,\boldsymbol{u}_{\ell,J},\ldots)\quad\\
\left.  \dfrac{\partial^{\ell}\boldsymbol{u}}{\partial t^{\ell}}\right\vert
_{t=t_{0}}=\boldsymbol{h}_{\ell}(\boldsymbol{x})
\end{array}
\right.  \quad\ell<k,\ |J|\ +\ell{}\leq k \label{Cauchy N}%
\end{equation}
where $\boldsymbol{f}$ and the $\boldsymbol{h}_{\ell}$'s are analytic vector
valued functions of their arguments. A system of $m$ PDEs in $m$ unknowns
functions is in \emph{normal form} if it is of the kind (\ref{N})

\emph{Cauchy-Kowalewski Theorem} (see, for instance, \cite{ch62}) asserts
that, locally, there exists a unique solution of the Cauchy problem (\ref{N})
+ (\ref{Cauchy N}). Since the initial data in (\ref{Cauchy N}) completely
determine the Taylor series of $\boldsymbol{u}$ at points of the \emph{initial
surface }$t=t_{0}$, the proof basically consists in checking convergence of
the series.

\subsubsection{General Cauchy problems}

Often, e.g. in relativistic theories, there is no preferred \textquotedblleft
space + time splitting\textquotedblright\ of the space-time. In this case, it
is generically advisable not to break the covariance by an arbitrary choice of
a time coordinate. Thus, a Cauchy problem is better posed on a generic
hypersurface in the space-time. Namely, consider a system of PDEs in $m$
unknown functions $\boldsymbol{u}$ in the general form
\begin{equation}
\boldsymbol{F}(x,\ldots,\boldsymbol{u}_{I},\ldots)=0,\quad|I|{}\leq k\label{1}%
\end{equation}
and a generic (Cauchy, i.e., initial) hypersurface
\[
\Sigma:z(x)=0,
\]
where $\boldsymbol{F}$ are smooth functions with independent differentials and
$z$ is a smooth function with non-vanishing gradient. The first normal
derivative of $\boldsymbol{u}$ at a point $(x^{1},\ldots,x^{n})$ of $\Sigma$
is
\[
\dfrac{\partial\boldsymbol{u}}{\partial z}:=\dfrac{\partial z}{\partial x^{1}%
}\boldsymbol{u}_{1}+\cdots+\dfrac{\partial z}{\partial x^{n}}\boldsymbol{u}%
_{n}.
\]
Put
\[
\dfrac{\partial^{\ell}\boldsymbol{u}}{\partial z^{\ell}}:=\dfrac{\partial
}{\partial z}\cdots\dfrac{\partial}{\partial z}\boldsymbol{u}.
\]
An initial data problem (Cauchy problem) on $\Sigma$ can be posed as follows:
\begin{equation}
\left\{
\begin{array}
[c]{l}%
\boldsymbol{F}(x,\ldots,\boldsymbol{u}_{I},\ldots)=0\quad\\
\left.  \dfrac{\partial^{\ell}\boldsymbol{u}}{\partial z^{\ell}}\right\vert
_{z=0}=\boldsymbol{h}_{\ell}(\boldsymbol{x})
\end{array}
\right.  \quad|I|{}\leq k,\quad\ell<k.\label{Cauchy}%
\end{equation}
If Problem (\ref{Cauchy}) could be recast in the normal form (\ref{Cauchy N}),
then, under additional analiticity condition, I could apply the
Cauchy-Kowalewski Theorem and get existsence and uniqueness of solutions. For
simplicity, I assume, from now on, that system (\ref{1}) is

\begin{enumerate}
\item \emph{weakly determined}, in the sense that it consists of precisely $m$ equations.

\item \emph{quasi-linear}, i.e.,
\begin{equation}
\boldsymbol{F}=\boldsymbol{A}^{j_{1}\cdots j_{k}}\cdot\boldsymbol{u}%
_{j_{1}\cdots j_{k}}+\boldsymbol{g},\label{6}%
\end{equation}

\end{enumerate}

where, for all multi-indexes $j_{1}\cdots j_{k}$,
\[
\boldsymbol{A}^{j_{1}\cdots j_{k}}=\boldsymbol{A}^{j_{1}\cdots j_{k}}%
(x,\ldots,\boldsymbol{u}_{J},\ldots),\quad|J|{}<k
\]
is an $m\times m$ matrix valued function, and\quad%
\[
\boldsymbol{g}=\boldsymbol{g}(x,\ldots,\boldsymbol{u}_{J},\ldots),\quad|J|{}<k
\]
is a vector valued function.

Notice that quasi-linearity is a condition invariant under a change of (both
independent and dependent) coordinates. In particular, it is easy to see (by
induction on $k$) that if $x=(x^{1},\ldots,x^{n})\longmapsto\bar{x}=(\bar
{x}^{1},\ldots,\bar{x}^{n})$ is a diffeomorphism, then
\[
\boldsymbol{F}=\boldsymbol{\bar{A}}{}^{j_{1}\cdots j_{k}}\cdot\boldsymbol{\bar
{u}}_{j_{1}\cdots j_{k}}+\boldsymbol{\bar{g}},
\]
where the $\boldsymbol{\bar{u}}_{I}$'s are derivatives of $\boldsymbol{u}$
with respect to the $\bar{x}$'s,
\[
\boldsymbol{\bar{A}}{}^{i_{1}\cdots i_{k}}=\dfrac{\partial\bar{x}^{i_{1}}%
}{\partial x^{j_{1}}}\cdots\dfrac{\partial\bar{x}^{i_{k}}}{\partial x^{j_{k}}%
}\boldsymbol{A}^{j_{1}\cdots j_{k}}\quad\text{and\quad}\boldsymbol{\bar{g}%
}=\boldsymbol{\bar{g}}(\bar{x},\ldots,\boldsymbol{\bar{u}}_{J},\ldots
),\quad|J|{}<k.
\]
In particular, the coefficients $\boldsymbol{A}^{i_{1}\cdots i_{k}}$ of the
highest order (linear) term transform as a contravariant symmetric tensor
under a change of (independent) coordinates. The contravariant tensor
\[
\boldsymbol{A}=(\boldsymbol{A}^{i_{1}\cdots i_{k}})
\]
is called the (principal) \emph{symbol} of the quasi-linear operator
$\boldsymbol{F}$ (see Subsection \ref{SecSymb} for an intrisic definition of
the symbol).

\begin{remark}
Limiting the discussion to determined, quasi-linear systems of PDEs is not
really restrictive for physical applications. Indeed, such systems are
particularly relevant in Physics, since Euler-Lagrange PDEs are precisely of
this form.
\end{remark}

Now, choose independent coordinates adapted to $\Sigma$, i.e., complete $z$ to
a system of coordinates $(z,y^{1},\ldots,y^{n-1})$. Then $\boldsymbol{y}%
=(y^{1},\ldots,y^{n-1})$ can be understood as internal coordinates on $\Sigma
$. In the new coordinates, Eq. (\ref{1}) becomes
\begin{equation}
\dfrac{\partial z}{\partial x^{j_{1}}}\cdots\dfrac{\partial z}{\partial
x^{j_{k}}}\boldsymbol{A}^{j_{1}\cdots j_{k}}\cdot\dfrac{\partial
^{k}\boldsymbol{u}}{\partial z^{k}}=\boldsymbol{{f}}(z,\boldsymbol{y}%
,\ldots,\boldsymbol{\bar{u}}_{\ell,J},\ldots),\quad\ell<k,\ |J|\ +\ell{}\leq
k.\label{Cauchy II}%
\end{equation}
for a suitable $\boldsymbol{f}$, where
\[
\boldsymbol{\bar{u}}_{\ell,J}=\frac{\partial^{|J|+\ell}\boldsymbol{u}%
}{\partial\boldsymbol{y}^{J}\partial z^{\ell}}.
\]
If
\begin{equation}
\left.  \det\left(  \dfrac{\partial z}{\partial x^{j_{1}}}\cdots
\dfrac{\partial z}{\partial x^{j_{k}}}\boldsymbol{A}^{j_{1}\cdots j_{k}%
}\right)  \right\vert _{z=0}\neq0\label{9}%
\end{equation}
then (\ref{Cauchy II}) can be clearly recast in the normal form
(\ref{Cauchy N}), around $\Sigma$. Notice, however, that the $\boldsymbol{A}%
^{j_{1}\cdots j_{k}}$ will generically depend on the $\boldsymbol{u}_{J}$,
$|J|{}<k$. Accordingly, \emph{unequality (\ref{9}) is actually a condition on
}$\Sigma$\emph{ and initial data on it, rather then a condition on the sole
}$\Sigma$.

\subsection{Characteristic Cauchy Data}

\subsubsection{Characteristic covectors and characteristic Cauchy data}

One is thus led to consider the $m\times m$ matrix
\[
\boldsymbol{A}(\boldsymbol{p})=p_{j_{1}}\cdots p_{j_{k}}\boldsymbol{A}%
^{j_{1}\cdots j_{k}}%
\]
for an arbitary co-vector $\boldsymbol{p}=p_{i}dx^{i}$. More precisely,
$\boldsymbol{A}(\boldsymbol{p})$ is a \emph{matrix-valued, homogeneous
polynomial function on cotangent spaces} to the space-time (see \cite{a04} for
a nice example). Notice that, in general, it does also depend on the
space-time point $x$ and on derivatives of the field at the point $x$ up to
the order $k-1$. For simplicity, I assume, temporarily, that $\boldsymbol{A}%
(\boldsymbol{p})$ is generically invertible, i.e.,
\[
\operatorname{rank}\boldsymbol{A}(\boldsymbol{p})=m
\]
somewhere, and therefore, almost everywhere, in the space of the
$\boldsymbol{p}$'s. Notice that $\det\boldsymbol{A}(\boldsymbol{p})$ is a
homogenous polynomial in the $p_{i}$'s. Therefore, if the equation
\[
\det\boldsymbol{A}(\boldsymbol{p})=0
\]
is compatible, then it determines a closed, nowhere dense, conic subset of the
space of the $\boldsymbol{p}$'s, called (up to projectivization) the
\emph{characteristic variety of the equation} (\ref{1}) (see, for instance,
\cite{b...91}). Points of the characteristic variety are called
\emph{characteristic covectors} and play an important role for different
aspects of the theory of PDEs, namely: the Cauchy problem and singularities of
solutions (as discussed below), the classification of PDEs \cite{v73}, the
method of intermediate integrals \cite{klv86} (for finding solutions of a PDE
by integrating lower order PDEs).

An hypersurface $\Sigma:z=0$ such that
\begin{equation}
\det\boldsymbol{A}(dz)|_{z=0}=\left.  \det\left(  \dfrac{\partial z}{\partial
x^{j_{1}}}\cdots\dfrac{\partial z}{\partial x^{j_{k}}}\boldsymbol{A}%
^{j_{1}\cdots j_{k}}\right)  \right\vert _{z=0}=0\label{CE}%
\end{equation}
is called a \emph{characteristic (Cauchy) surface} (see, for instance,
\cite{ch62}). Beware, however, that this is an abuse of terminology (mutuated
by the theory of \emph{linear PDEs}). Indeed, as already remarked, (\ref{CE})
is actually a condition on $\Sigma$ and initial data on it. Accordingly, we
should rather speak about \emph{characteristic (Cauchy) data}. The initial
value problem may not be well-posed (i.e., there may be no existence and
uniqueness, even for analytic data), in general, on characteristic surfaces.
In particular, initial data are \textquotedblleft\emph{constrained}%
\textquotedblright\ on a characteristic surface, in the sense that not all
initial data on a characteristic surface are \emph{admissible}, i.e., are
compatible with the PDE. To see this, let $\Sigma$ be characteristic, and
$q=m-\operatorname{rank}\boldsymbol{A}(dz)|_{z=0}>0$. As a minimal regularity
condition, I assume $q$ to be constant on $\Sigma$. Then there is a non zero
$[q,m]$ matrix $\boldsymbol{M}=\boldsymbol{M}(z,\boldsymbol{y},\ldots
,\boldsymbol{\bar{u}}_{\ell,J},\ldots)$ such that
\[
(\boldsymbol{M}\cdot\boldsymbol{A}(dz))|_{z=0}=0.
\]
It follows that Eq.{} (\ref{Cauchy II}) may only possess solution if
\begin{equation}
(\boldsymbol{M}\cdot\boldsymbol{f})(z,\boldsymbol{y},\ldots,\boldsymbol{\bar
{u}}_{\ell,J},\ldots)|_{z=0}=0.\label{Sigma}%
\end{equation}
This last equation may be interpreted as a system of (generically non-linear)
PDEs constraining the initial data
\[
\left.  \frac{\partial^{\ell}\boldsymbol{u}}{\partial z^{\ell}}\right\vert
_{z=0}%
\]
on the characteristic surface $\Sigma$.

Finally, notice that (\ref{CE}) may be interpreted as a first order,
polynomial, PDE whose unknown is a hypersurface in the space-time. To see
this, assume $z$ to be in the form $z=t-\tau(\boldsymbol{x})$, and
$\boldsymbol{y}=\boldsymbol{x}$. Then (\ref{CE}) becomes
\begin{equation}
\det\left(
%TCIMACRO{\tsum \nolimits_{\ell=0}^{k}}%
%BeginExpansion
{\textstyle\sum\nolimits_{\ell=0}^{k}}
%EndExpansion
\frac{\partial\tau}{\partial\boldsymbol{y}^{a_{1}}}\cdots\frac{\partial\tau
}{\partial\boldsymbol{y}^{a_{\ell}}}\boldsymbol{B}^{a_{1}\cdots a_{\ell}%
}\right)  =0\label{CEdiv}%
\end{equation}
where
\[
\boldsymbol{B}^{a_{1}\cdots a_{\ell}}=\frac{k!}{\ell!}\boldsymbol{A}%
^{a_{1}\cdots a_{\ell}n\cdots n},\quad a_{1},\ldots,a_{\ell}=1,\ldots,n-1.
\]
It should be stressed, however, that the \textquotedblleft
coefficients\textquotedblright\ $\boldsymbol{B}^{a_{1}\cdots a_{\ell}}$ depend
in general on $\boldsymbol{u}$ and its derivatives up to the order $k-1$. When
the $\boldsymbol{B}^{a_{1}\cdots a_{\ell}}$'s do only depend on independent
variables $\boldsymbol{y}$, (e.g., when Eq.{} (\ref{1}) is linear) Eq.{}
(\ref{CEdiv}) is a first order, scalar (inhomogeneous polynomial) PDE in the
unknown $\tau$ that can be treated, for instance, with the method of
characteristics. In this case, one usually refers to \emph{characteristic
lines} of (\ref{CEdiv}) as \emph{bicharacteristics} of (\ref{1}). We will come
back to bicharacteristics (and the method of characteristics for scalar PDEs)
in Section \ref{SecBic}.

\subsubsection{Physical interpretation of characteristic
surfaces\label{SecPhys}}

On another hand, \emph{singularities of solutions of a system of PDEs occur
along characteristic surfaces}. To clarify this sentence, let $\boldsymbol{u}%
_{0}$ be a fiducial, background solution of (\ref{1}), and $\Sigma$ be an
hypersurface bounding a region $\Omega_{0}$ of the space-time. We search for a
(possibly singular along $\Sigma$) solution $\boldsymbol{u}$ of (\ref{1})
which agrees with $\boldsymbol{u}_{0}$ in $\Omega_{0}$ but is everywhere
different from $\boldsymbol{u}_{0}$ outside $\Sigma$. On \emph{physical
ground}, I assume that all derivatives of $\boldsymbol{u}$ up to the order
$k-1$ are continuous along $\Sigma$. In particular, $\boldsymbol{u}$ and
$\boldsymbol{u}_{0}$ are both solutions of a Cauchy problem of the form
(\ref{Cauchy}). It follows that $\Sigma$ must be a characteristic surface. In
other words, as already stated, singularities (e.g., wave fronts) occur along
characteristic surfaces, i.e., the \emph{boundary of a disturbance in the
space-time is a characteristic surface }\cite{lc31}. Thus, to understand how
disturbaces of a specific field (with specific field equations) propagate in
the space time, one has to solve the characteristic equation (\ref{CE}).

Under suitable conditions, a characteristic surface is actually equipped with
a field of directions that integrates to a 1 dimensional foliation, whose
leaves are traditionally referred to as \emph{bicharacteristics} (Hadamard
terminology). Accordingly, \emph{singularities of solutions propagate along
bicharacteristics}. From a physical point of view, one may interpret
characteristic surfaces as \emph{wave-fronts} and bicharacteristics as
\emph{rays}. Under this interpretations the passage from characteristics to
bicharacteristics is the \emph{passage from wave-optics to geometric optics}
(see \cite{l64} for more details, see also \cite{gs77}). Alternatively,
bicharacteristics can be interpreted as trajectories of particles. If one
adopts this interpretation, they describe the motion of a particle-like
counterpart of the field under consideration. This relates \emph{the principle
of wave-particle duality} to the \emph{geometric theory of PDEs}.

\subsubsection{Characteristics of Euler-Lagrange equations}

In general, $\boldsymbol{A}(\boldsymbol{p})$ may be non-invertible everywhere
on the space of $\boldsymbol{p}$'s. In this case, initial data are constrained
on every Cauchy hypersurface and the Cauchy-Kowalewski theorem fails. Let $r$
be the maximum rank of $\boldsymbol{A}(\boldsymbol{p})$ on the space of
$\boldsymbol{p}$'s. Then $\operatorname{rank}\boldsymbol{A}(\boldsymbol{p})=r$
almost everywhere in the space of $\boldsymbol{p}$'s. In this case, one define
a characteristic surface to be an hypersurface $\Sigma:z=0$ such that
\[
\operatorname{rank}\boldsymbol{A}(dz)|_{z=0}=\left.  \operatorname{rank}%
\left(  \dfrac{\partial z}{\partial x^{j_{1}}}\cdots\dfrac{\partial
z}{\partial x^{j_{k}}}\boldsymbol{A}^{j_{1}\cdots j_{k}}\right)  \right\vert
_{z=0}<r.
\]
In this general case, characteristic surfaces still play a role in Cauchy
problems and the theory of sigularity propagation, but I will not enter this
here. However, notice that the general situation does occur for the field
equation of a gauge theory as I briefly discuss now. Indeed, let (\ref{1}) be
the Euler-Lagrange (EL) equations determined by a variational principle
\[
\int L(x^{1},\ldots,x^{n},\ldots,\boldsymbol{u}_{J},\ldots)d^{n}x,\quad
|J|{}\leq\ell
\]
where $L(x^{1},\ldots,x^{n},\ldots,\boldsymbol{u}_{J},\ldots)d^{n}x$ is a
Lagrangian density depending on derivatives of the fields up to the order
$\ell$. Then
\[
\boldsymbol{F}=\sum_{|J|{}\leq\ell}(-)^{|J|}D_{J}\dfrac{\partial L}%
{\partial\boldsymbol{u}_{J}},
\]
where, for $J=j_{1}\cdots j_{r}$, $D_{J}:=D_{j_{1}}\circ\cdots\circ D_{j_{r}}%
$, and
\[
D_{i}:=\dfrac{\partial}{\partial x^{i}}+\sum_{I}\boldsymbol{u}_{Ii}\cdot
\dfrac{\partial}{\partial\boldsymbol{u}_{I}}%
\]
is the $i$-th total derivative. Then
\[
\boldsymbol{F}=\sum\limits_{j_{1}\leq\cdots\leq j_{\ell}}\sum\limits_{k_{1}%
\leq\cdots\leq k_{\ell}}\dfrac{\partial^{2}L}{\partial\boldsymbol{u}%
_{j_{1}\cdots j_{\ell}}\partial\boldsymbol{u}_{k_{1}\cdots k_{\ell}}}%
\cdot\boldsymbol{u}_{j_{1}\cdots j_{\ell}k_{1}\cdots k_{\ell}}+\boldsymbol{g}%
\]
where\quad%
\[
\boldsymbol{g}=\boldsymbol{g}(x^{1},\ldots,x^{n},\ldots,\boldsymbol{u}%
_{J},\ldots),\quad|J|{}<2\ell.
\]
For a gauge invariant Lagrangian
\[
\det\boldsymbol{A}(\boldsymbol{p})=\det\left(  \sum\limits_{j_{1}\leq
\cdots\leq j_{\ell}}\sum\limits_{k_{1}\leq\cdots\leq k_{\ell}}\dfrac
{\partial^{2}L}{\partial\boldsymbol{u}_{j_{1}\cdots j_{\ell}}\partial
\boldsymbol{u}_{k_{1}\cdots k_{\ell}}}p_{j_{1}}\cdots p_{j_{\ell}}p_{k_{1}%
}\cdots p_{k_{\ell}}\right)  =0
\]
for all $\boldsymbol{p}$'s (see examples below).

\subsubsection{Characteristics of fully non-linear equations\label{SecNL}}

Finally, I briefly discuss the case when Eq.{} (\ref{1}) is not quasi-linear.
In this case, a careful use of the inverse function theorem shows that the
Cauchy problem is well posed on any hypersurface $\Sigma:z=0$ such that
\[
\left.  \det\left(  \sum\nolimits_{j_{1}\leq\cdots\leq j_{k}}\dfrac{\partial
z}{\partial x^{j_{1}}}\cdots\dfrac{\partial z}{\partial x^{j_{k}}}%
\frac{\partial\boldsymbol{F}}{\partial\boldsymbol{u}_{j_{1}\cdots j_{k}}%
}\right)  \right\vert _{z=0}\neq0,
\]
Accordingly, all the above considerations remain valid up to a substitution
\[
\boldsymbol{A}^{j_{1}\cdots j_{k}}\longrightarrow\frac{(j_{1}\cdots j_{k}%
)!}{k!}\,\frac{\partial\boldsymbol{F}}{\partial\boldsymbol{u}_{j_{1}\cdots
j_{k}}},
\]
where $(i_{1}\cdots i_{k})!/k!$ is a suitable combinatorial coefficients that
accounts for the fact that
\[
\boldsymbol{u}_{j_{1}\cdots j_{\ell}}=\boldsymbol{u}_{j_{\sigma(1)}\cdots
j_{\sigma(\ell)}}%
\]
for every permutation $\sigma$ of $\{1,\ldots,\ell\}$. Specifically, let
$j=1,\ldots,n$ appear $N_{j}$ times in the multi-index $j_{1}\cdots j_{k}$.
Then $(j_{1}\cdots j_{k})!:=N_{1}!\cdots N_{n}!$. If the matrix
\[
\sum\nolimits_{j_{1}\leq\cdots\leq j_{k}}p_{j_{1}}\cdots p_{j_{k}}%
\frac{\partial\boldsymbol{F}}{\partial u_{j_{1}\cdots j_{k}}}%
\]
is generically invertible on the space of $\boldsymbol{p}$'s, then
$\Sigma:z=0$ is a characteristic surface if
\[
\left.  \det\left(  \sum\nolimits_{j_{1}\leq\cdots\leq j_{k}}\dfrac{\partial
z}{\partial x^{j_{1}}}\cdots\dfrac{\partial z}{\partial x^{j_{k}}}%
\frac{\partial\boldsymbol{F}}{\partial u_{j_{1}\cdots j_{k}}}\right)
\right\vert _{z=0}=0.
\]

\subsection{Examples}

\subsubsection{Klein-Gordon and wave equations on a curved space-time}

Let $\boldsymbol{g}=g_{ij}dx^{i}dx^{j}$ be a Riemannian, or pseudo-Riemannian
metric on an open subset $U$ of $\mathbb{R}^{n}$. Consider the following
linear equation
\begin{equation}
g^{ij}\nabla_{i}\nabla_{j}u=0,\label{2}%
\end{equation}
where $u$ is an unknown function on $U$, and $\nabla$ is the Levi-Civita
connection of $\boldsymbol{g}$. Eq. (\ref{2}) is the EL equation coming from
the action functional%
\[
-\frac{1}{2}\int(g^{ij}\nabla_{i}u\nabla_{j}u)\sqrt{|\det\boldsymbol{g}|}%
d^{n}x.
\]
The symbol of the\emph{ }operator $\boldsymbol{F}=g^{ij}\nabla_{i}\nabla_{j}$
is
\[
\boldsymbol{A}=\boldsymbol{g}^{-1}=:(g^{ij}).
\]
Accordingly,
\[
\boldsymbol{A}(\boldsymbol{p})=g^{ij}p_{i}p_{j}=\boldsymbol{g}^{-1}%
(\boldsymbol{p},\boldsymbol{p}),
\]
which is generically non-zero. The characteristic variety is the quadric
\[
\boldsymbol{A}(\boldsymbol{p})=g^{ij}p_{i}p_{j}=\boldsymbol{g}^{-1}%
(\boldsymbol{p},\boldsymbol{p})=0,
\]
and characterictic surfaces $\Sigma:z=0$ are defined by
\[
\left.  g^{ij}\frac{\partial z}{\partial x^{i}}\frac{\partial z}{\partial
x^{j}}\right\vert _{z=0}=\boldsymbol{g}^{-1}(\boldsymbol{d}z,\boldsymbol{d}%
z)|_{z=0}=0.
\]
In particular, if $\boldsymbol{g}$ is Riemannian then $\boldsymbol{F}=\Delta$
is the (curved) Laplacian, which is an elliptic operator, $\boldsymbol{F}=0$
is the (curved) Laplace equation, and there are no characteristic surfaces. On
another hand, if $\boldsymbol{g}$ is Lorentzian, $\boldsymbol{F}=\square$ is
the (curved) d'Alambertian, which is a hyperbolic operator, $\boldsymbol{F}=0$
is the (curved) wave equation, and characteristic surfaces are precisely the
null hypersurfaces. In this case, one concludes that wave fronts are
light-like hypersurfaces. Notice that the (curved) Klein-Gordon operator
$\square+m^{2}$ has the same symbol as the d'Alambertian, and, therefore, the
Klein-Gordon equation has the same characteristic surfaces as the wave equation.

\subsubsection{Dirac equation on Minkowski space-time}

Let $\boldsymbol{\eta}=\eta_{\mu\nu}dx^{\mu}dx^{\nu}$ be the Minkowski metric
on $\mathbb{R}^{4}$. The Dirac equation is the linear, first order (system of)
PDE(s) given by
\[
(i\boldsymbol{\gamma}^{\mu}\partial_{\mu}-m)\boldsymbol{u}=0,
\]
where $\boldsymbol{u}=(u^{0},u^{1},u^{2},u^{3})$ is a 4-component (complex)
spinor and $\boldsymbol{\gamma}^{0},\boldsymbol{\gamma}^{1},\boldsymbol{\gamma
}^{2},\boldsymbol{\gamma}^{3}$ are the $4\times4$ Dirac matrices. The symbol
of the operator $\boldsymbol{F}=i\boldsymbol{\gamma}^{\mu}\partial_{\mu}-m$
is
\[
\boldsymbol{A}=(i\boldsymbol{\gamma}^{\mu}).
\]
Accordingly,
\[
\boldsymbol{A}(\boldsymbol{p})=i\boldsymbol{\gamma}^{\mu}p_{\mu}%
\]
which is generically invertible. The characteristic variety is defined by the
4-th order algebraic equation
\[
\det\boldsymbol{A}(\boldsymbol{p})=\det(i\boldsymbol{\gamma}^{\mu}p_{\mu})=0.
\]
An easy computation (first performed by G. Racah in the 30th's \cite{r32})
shows that
\[
\det\boldsymbol{A}(\boldsymbol{p})=(\eta^{\mu\nu}p_{\mu}p_{\nu})^{2}%
=(\boldsymbol{\eta}^{-1}(\boldsymbol{p},\boldsymbol{p}))^{2}=0,
\]
and characteristic surfaces $\Sigma:z=0$ are defined by
\[
\left.  \left(  \eta^{\mu\nu}\frac{\partial z}{\partial x^{\mu}}\frac{\partial
z}{\partial x^{\nu}}\right)  ^{2}\right\vert _{z=0}=(\boldsymbol{\eta}%
^{-1}(\boldsymbol{d}z,\boldsymbol{d}z))^{2}|_{z=0}=0.
\]
Therefore, characteristic surfaces of the Dirac equations are precisely null
surfaces in the Minkowski space-time (and Racah himself interpreted this
result in terms of the Heisenberg principle).

\subsubsection{Maxwell equations on a curved space-time}

Let $\boldsymbol{g}=g_{ij}dx^{i}dx^{j}$ be a Lorentzian metric on an open
subset $U$ of $\mathbb{R}^{4}$. The (vacuum) Maxwell equations in $U$ read
\[
F_{j}=g^{ik}\nabla_{k}(\nabla_{i}u_{j}-\nabla_{j}u_{i})=0,
\]
where $\boldsymbol{u}=(u_{0},u_{1},u_{2},u_{3})$ are the components of a
differential $1$-form (the electromagnetic potential) on $U$. Maxwell
equations are the EL equations coming from the action functional
\[
-\int g^{ik}g^{j\ell}\nabla_{\lbrack i}u_{j]}\nabla_{\lbrack k}u_{\ell]}%
\sqrt{|\det\boldsymbol{g}|}d^{4}x.
\]
Now
\[
F_{j}=\left(  g^{ik}\delta_{j}^{\ell}-g^{\ell k}\delta_{j}^{i}\right)
\frac{\partial^{2}u_{\ell}}{\partial x^{k}\partial x^{i}}+\cdots,
\]
where the dots $\cdots$ denote lower order terms. Accordingly,
\[
\boldsymbol{A}(\boldsymbol{p})_{j}^{\ell}=\left(  g^{ik}\delta_{j}^{\ell
}-g^{\ell k}\delta_{j}^{i}\right)  p_{k}p_{i},
\]
i.e.,
\[
\boldsymbol{A}(\boldsymbol{p})=\boldsymbol{g}^{-1}(\boldsymbol{p}%
,\boldsymbol{p)}\,\mathbf{I{}-{}}\boldsymbol{p}^{\sharp}\otimes\boldsymbol{p},
\]
where $\boldsymbol{p}^{\sharp}:=\boldsymbol{g}^{-1}(\boldsymbol{p}%
,-\boldsymbol{)}$. Notice that $\boldsymbol{A}(\boldsymbol{p})$ is never
invertible. Indeed, $\operatorname{rank}\boldsymbol{A}(\boldsymbol{p})$ is
generically $3$ rather then $4$. This corresponds to the fact that gauge
freedom is parametrized by $1$ arbitrary function on the space-time. In this
(degenerate) case, characteristic surfaces $\Sigma:z=0$ are defined by
\[
\operatorname{rank}\boldsymbol{A}(dz)|_{z=0}<3.
\]
But $\operatorname{rank}\boldsymbol{A}(\boldsymbol{p})<3$ iff $\boldsymbol{p}$
is a null covector, and, in this case, $\operatorname{rank}\boldsymbol{A}%
(\boldsymbol{p})=1$ whenever $\boldsymbol{p}\neq0$. One concludes that the
characteristic surfaces of Maxwell equations in curved space-time are again
null hypersurfaces \cite{w28}.

Notice that the degeneracy of the matrix $\boldsymbol{A}(\boldsymbol{p})$ can
be cured by gauge fixing. For instance, for Maxwell equations in the Lorentz
gauge ($\nabla^{i}u_{i}=0$)
\[
\boldsymbol{A}(\boldsymbol{p})=\boldsymbol{g}^{-1}(\boldsymbol{p}%
,\boldsymbol{p)}\,\mathbf{I{}},
\]
which is generically invertible, and degenerates iff $\boldsymbol{g}%
^{-1}(\boldsymbol{p},\boldsymbol{p)}=0$ again \cite{w28}.

\subsubsection{Einstein equations}

Let $U$ be an open subset of $\mathbb{R}^{4}$. The (vacuum) Einstein equations
in $U$ read
\[
\mathbf{Ric}[\boldsymbol{u}]=0,
\]
where $\boldsymbol{u}=u_{ij}dx^{i}dx^{j}$ is an unknown Lorentzian metric on
$U$ and $\mathbf{Ric}[\boldsymbol{u}]=R_{ij}[\boldsymbol{u}]dx^{i}dx^{j}$ is
its Ricci tensor. Einstein equations are the EL equations coming from the
action functional
\[
\int u^{ij}R_{ij}[\boldsymbol{u}]\sqrt{|\det\boldsymbol{u}|}d^{4}x.
\]
The symbol of the Ricci operator $\mathbf{Ric}$ has been first computed by
Levi-Civita in the 30th's. One has
\[
R_{ij}[\boldsymbol{u}]=\left(  2\delta{}_{i}^{[m}u_{{}}^{k][\ell}\delta{}%
_{j}^{n]}-u_{{}}^{n[m}u_{{}}^{k]\ell}u_{ij}\right)  \frac{\partial^{2}%
u_{k\ell}}{\partial x^{m}\partial x^{n}}+\cdots
\]
where the dots $\cdots$ denote lower order terms. Accordingly,
\[
\boldsymbol{A}(\boldsymbol{p})_{ij}^{k\ell}=\left(  2\delta{}_{i}^{[m}u_{{}%
}^{k][\ell}\delta{}_{j}^{n]}-u_{{}}^{n[m}u_{{}}^{k]\ell}g_{ij}\right)
p_{k}p_{\ell},
\]
which should be understood as entries of a $10\times10$ matrix (the pairs $ij$
and $k\ell$ are to be ordered, for instance, lexicographically). Levi-Civita
proved that $\boldsymbol{A}(\boldsymbol{p})$ is never invertible, and
$\operatorname{rank}\boldsymbol{A}(\boldsymbol{p})$ is generically $6$ rather
then $4$. This corresponds to the fact that the gauge freedom is parametrized
by $4$ arbitrary functions on the space-time. Finally, $\operatorname{rank}%
\boldsymbol{A}(\boldsymbol{p})<6$ iff
\[
\boldsymbol{u}^{-1}(\boldsymbol{p},\boldsymbol{p})=u^{ij}p_{i}p_{j}=0,
\]
and, in this case, $\operatorname{rank}\boldsymbol{A}(\boldsymbol{p})=4$
whenever $\boldsymbol{p}\neq0$. One concludes that the characteristic surfaces
of Einstein equations are null hypersurfaces with respect to the unknown
metric $\boldsymbol{u}$. This is a typical example when (\ref{CE}) is a
condition on $\Sigma$ and initial data on it (in this case, the metric on it)
and not only on $\Sigma$ itself. Notice that, from a physical point of view,
the outcome of this and the previous three subsections is that \emph{the phase
velocity of gravitational, electromagnetic, Dirac, and Klein-Gordon field is
the speed of light.}

\subsubsection{An unphysical, fully non-linear example}

Consider the scalar PDE in two independent variables $x,y$:
\begin{equation}
u_{yyy}-\left(  u_{xxy}\right)  ^{2}+u_{xxx}u_{xyy}=0.\label{15}%
\end{equation}
Eq.{} (\ref{15}) is a third order \emph{Monge-Amp\`{e}re equation \cite{b92}}.
For $\boldsymbol{p}=pdx+qdy$, one has
\[
\boldsymbol{A}(\boldsymbol{p})=u_{xyy}p^{3}-2u_{xxy}p^{2}q+u_{xxx}pq^{2}%
+q^{3},
\]
which is generically non-zero. Accordingly, a hypersurface $\Sigma:z=z(x,y)$
is characteristic iff
\[
u_{xyy}z_{x}^{3}-2u_{xxy}z_{x}^{2}z_{y}+u_{xxx}z_{x}z_{y}^{2}+z_{y}^{3}=0.
\]
Notice that $z_{x}\neq0$, otherwise $z_{x}=z_{y}=0$. Therefore, one can search
for $\Sigma$ in the form $\Sigma:x=\tau(y)$, which gives
\[
u_{xyy}+2u_{xxy}\tau_{y}+u_{xxx}\tau_{y}^{2}-\tau_{y}^{3}=0
\]
depending on (constrained) initial data on $\Sigma$.

\section{ Singularities of Solutions of PDEs}

\subsection{PDEs and Jet Spaces}

Most of the considerations done in the previous lecture are independent of the
choice of coordinates. This suggests that there is an intrinsic, geometric
theory capturing the concept of characteristics of a system of PDEs. This is
actually the case. The aim of this section is to provide a gentle introduction
to basics of the geometric theory of (nonlinear) PDEs, their characteristics,
and (fold-type) singularities of their solutions. In particular, I will
present a rigorous, mathematical version of the physical considerations in
Subsection \ref{SecPhys}. The main results will be presented without a proof
and the interested reader should refer to the bibliography for details.
Indeed, a deeper analysis would show that many branches of Mathematics enter
the intrinsic theory of PDEs, namely: differential geometry and differential
topology, commutative algebra and algebraic geometry, homological algebra and
algebraic topology.

I begin with a geometric framework for PDEs, namely, jet spaces (for more
details about jet spaces, see \cite{b...99,klv86}).

\subsubsection{Jets of sections}

Let $\pi:E\longrightarrow M$ be a fiber bundle, and let $(x^{1},\ldots
,x^{n},\boldsymbol{u})$ be a bundle chart on $E$, i.e., $(x^{1},\ldots,x^{n})$
are coordinates on $M$, and $\boldsymbol{u}=(u^{1},\ldots,u^{m})$ are fiber
coordinates on $E$. the $x$'s will be interpreted as independent variables,
and $\boldsymbol{u}$ as a set of dependent variables. From a physical point of
view, $M$ will be often interpreted as the space-time and sections of $\pi$ as
configurations of a field on it. I want to discuss PDEs imposed on sections of
$\pi$. To do this in a way which is manifestly independent of coordinates (and
any other auxiliary structure on $\pi$, e.g., a connection) it is necessary to
introduce jet spaces.

Two local sections $\sigma_{1}$ and $\sigma_{2}$ of $\pi$, locally given by
\[
\sigma_{1,2}:\boldsymbol{u}=\boldsymbol{f}_{1,2}(x^{1},\ldots,x^{n}),
\]
are \emph{tangent up to the order }$k$ at a point $x_{0}\in M$ if the $k$-th
order Taylor polynomials of $\boldsymbol{f}_{1}$ and $\boldsymbol{f}_{2}$
coincide at $x_{0}\equiv(x_{0}^{1},\ldots,x_{0}^{n})$:
\[
\frac{\partial^{|I|}\boldsymbol{f}_{1}}{\partial x^{I}}(x_{0}^{1},\ldots
,x_{0}^{n})=\frac{\partial^{|I|}\boldsymbol{f}_{2}}{\partial x^{I}}(x_{0}%
^{1},\ldots,x_{0}^{n}),\quad|I|{}\leq k.
\]
Tangency up to the order $k$ at $x_{0}$ is a well defined equivalence
relation. In particular, it is independent of coordinates. Denote by
$J_{x_{0}}^{k}\pi$ the set of equivalence classes. Finally put
\[
J^{k}\pi:=\coprod\limits_{x_{0}\in M}J_{x_{0}}^{k}\pi.
\]
It is called the $k$\emph{-th jet space of the bundle }$\pi$ and can be given
a canonical structure of smooth manifold as follows. First of all, for a local
section $\sigma$ of $\pi$, denote by $[\sigma]_{x}^{k}$ its class of tangency
up to the order $k$ at the point $x\in M$. It is a point of $J^{k}\pi$, which
is called the $k$\emph{-th jet of }$\sigma$\emph{ at }$x$, and can be
intepreted as (an intrinsic version of) the $k$-th order Taylor polynomial of
$\sigma$ at $x$. Notice that $J^{0}\pi$ identifies canonically with $E$.
Moreover, there are canonical surjections
\[
\pi_{k,\ell}:J^{k}\pi\longrightarrow J^{\ell}\pi,\quad k\geq\ell
\]
which consist in forgetting derivatives of order higher than $\ell$. There are
also surjections
\[
\pi_{k}:J^{k}\pi\longrightarrow M,\quad\lbrack\sigma]_{x}^{k}\longmapsto x.
\]
Let $U$ be a bundle coordinate domain in $E=J^{0}\pi$. On $\pi_{k,0}^{-1}(U)$
there are coordinates $(x^{1},\ldots,x^{n},\ldots,\boldsymbol{u}_{I},\ldots)$
given by
\[
\boldsymbol{u}_{I}([\sigma]_{x}^{k}):=\frac{\partial^{|I|}\boldsymbol{f}%
}{\partial x^{I}}(x^{1},\ldots,x^{n}),\quad|I|{}\leq k,
\]
where $\sigma$ is a local section of $\pi$ which in coordinates look as
\begin{equation}
\sigma:\boldsymbol{u}=\boldsymbol{f}(x^{1},\ldots,x^{n}).\label{3}%
\end{equation}
It is easy to see that the $J^{k}\pi$, with these coordinates, are smooth
manifolds, and the $\pi_{k,\ell}$'s (and, consequently, the $\pi_{k}$'s) are
fiber bundles.

A section $\sigma$ of $\pi$ can be \textquotedblleft
prolonged\textquotedblright\ to a section $j^{k}\sigma$ of $\pi_{k}$, its
$k$\emph{-th jet prolongation}, by putting
\[
(j^{k}\sigma)(x):=[\sigma]_{x}^{k}.
\]
If $\sigma$ looks locally as (\ref{3}), then $j^{k}\sigma$ looks locally as
\[
j^{k}\sigma:\boldsymbol{u}_{I}:=\frac{\partial^{|I|}\boldsymbol{f}}{\partial
x^{I}}(x^{1},\ldots,x^{n}),\quad|I|{}\leq k.
\]
Thus $j^{k}\sigma$ is a coordinate free version of \textquotedblleft partial
derivative functions of $\sigma$ up to the order $k$\textquotedblright. Notice
that not all sections of $\pi_{k}$ are of the form $j^{k}\sigma$. The latter
are sometimes called \emph{holonomic sections}.

\subsubsection{The Cartan Distribution}

In the following, I will denote simply by $J^{k}$ the space of $k$-th jets of
sections of $\pi$, if there is no risk of confusion. There is a canonical
structure on $J^{k}$, namely a distribution, which, in a sense, encodes the
\textquotedblleft differential relations among the jet coordinates
$\boldsymbol{u}_{I}$\textquotedblright. Let us fix a point $\theta\in J^{k}$.
If $\theta$ is the $k$-th jet at a point $x\in M$ of a section $\sigma$ of
$\pi$, then, clearly, $\operatorname{im}j^{k}\sigma\ni\theta$. Consider the
tangent space $R[\sigma]\subset T_{\theta}J^{k}$ to $\operatorname{im}%
j^{k}\sigma$ at $\theta$. Any subspace of $T_{\theta}J^{k}$ of the form
$R[\sigma]$ is called an $R$\emph{-plane} at $\theta$. Notice that $R[\sigma]$
does only depend on the $(k+1)$-th jet $\theta^{\prime}=[\sigma]_{x}^{k+1}$ of
$\sigma$. Accordingly, it will be denoted by $R_{\theta^{\prime}}$. The
correspondence
\[
\pi_{k+1,k}^{-1}(\theta)\longrightarrow\{R\text{-planes at }\theta
\},\quad\theta^{\prime}\longmapsto R_{\theta^{\prime}},
\]
is a bijection that allows to construct jet spaces inductively from lower
order ones.

$R$-planes at $\theta$ span a distinguished subspace $\mathcal{C}_{\theta}$ in
$T_{\theta}J^{k}$ and the correspondence $\mathcal{C}:\theta\longmapsto
\mathcal{C}_{\theta}$ is a smooth distribution on $J^{k}$ often called the
\emph{Cartan distribution}. The Cartan distribution is locally spanned by
vector fields
\[
D_{i}:=\dfrac{\partial}{\partial x^{i}}+\sum_{|I|{}<k}\boldsymbol{u}%
_{Ii}\dfrac{\partial}{\partial\boldsymbol{u}_{I}},\quad\dfrac{\partial
}{\partial\boldsymbol{u}_{J}},\quad i=1,\ldots,n,\quad|J|{}=k.
\]
For obvious reasons, the $D_{i}$'s are called the \emph{total derivatives}.
Dually, $\mathcal{C}$ is annihilated by $1$-forms
\[
\boldsymbol{\omega}_{I}:=d\boldsymbol{u}_{I}-\boldsymbol{u}_{Ii}dx^{i}%
,\quad|I|{}<k,
\]
called the \emph{Cartan forms. }This shows that $\mathcal{C}$\emph{ is not
involutive} (and, therefore, \emph{not integrable}). Actually, in general, the
Cartan distribution possesses many different (locally) maximal integral
submanifolds (even of different dimensions) through any point. For instance,
fibers of $\pi_{k,k-1}$ and images of holonomic sections are both maximal
integral submanifolds and there are more maximal integral submanifolds of
different kinds. However, if a maximal integral submanifold is horizontal with
respect to the projection $\pi_{k,k-1}$, then it is the image of a holonomic
section. In this sense the Cartan distribution \emph{\textquotedblleft
detects\textquotedblright\ holonomic sections}.

\subsubsection{Differential Equations}

Jet spaces formalize in a coordinate free way the concept of partial
derivatives. Accordingly, they allow a coordinate free definition of system of
PDEs. Specifically, a \emph{system of (non-linear) PDEs of the order }%
$k$\emph{ imposed on sections of the bundle }$\pi$ (in the following, simply a
PDE) is a submanifold $\mathcal{E}$ of $J^{k}$. Indeed $\mathcal{E}$ looks
locally as
\begin{equation}
\mathcal{E}:\boldsymbol{F}(x^{1},\ldots,x^{n},\ldots,\boldsymbol{u}_{I}%
,\ldots)=0,\quad|I|{}\leq k,\label{5}%
\end{equation}
which is a system of PDEs in the analytic sense up to the interpretation of
the $\boldsymbol{u}_{I}$'s as partial derivatives of the $\boldsymbol{u}$'s.
In view of (\ref{5}), it is meaningful to say that a a system of PDEs
$\mathcal{E}\subset J^{k}$ is \emph{(weakly) determined} (i.e., the number of
equations coincides with the number of dependent variables) if
$\operatorname{codim}\mathcal{E}=m$. The coordinate free definition of
\emph{solutions of }$\mathcal{E}$ should be now clear: a \emph{solution of
}$\mathcal{E}$ is a section $\sigma$ of $\pi$ such that $j^{k}\sigma$ takes
values in $\mathcal{E}$. Indeed, if $\sigma$ is locally given by (\ref{3}),
then the condition $\operatorname{im}j^{k}\sigma\subset\mathcal{E}$ reads
locally
\[
\boldsymbol{F}\left(  x^{1},\ldots,x^{n},\ldots,\partial^{|I|}\boldsymbol{f/}%
\partial x^{I},\ldots\right)  =0,\quad|I|{}\leq k
\]
which is a system of PDEs imposed on the $\boldsymbol{f}$'s. On a PDE
$\mathcal{E}\subset J^{k}$ one can consider the distribution $\mathcal{C}%
(\mathcal{E}):\theta\longmapsto\mathcal{C}_{\theta}(\mathcal{E}):=\mathcal{C}%
_{\theta}\cap T_{\theta}\mathcal{E}$. Under suitable regularity conditions on
$\mathcal{E}$, $\mathcal{C}(\mathcal{E})$ is a smooth distribution. It is then
clear that if a maximal integral submanifold of $\mathcal{C}(\mathcal{E})$ is
horizontal with respect to $\pi_{k,k-1}$, then it is the image of $j^{k}%
\sigma$ for some \emph{solution} $\sigma$ of $\mathcal{E}$. In other words,
\emph{smooth solutions of }$\mathcal{E}$\emph{ are in one-to-one
correspondence with maximal integral submanifolds of }$\mathcal{C}%
(\mathcal{E})\ $\emph{satisfying a suitable horizontality conditions}. The
main point here is that, \emph{relaxing this horizontality condition, one can
describe, in purely geometric terms, solutions with (specific type of)
singularities}.

\subsubsection{Jets of Submanifolds}

Notice that, in differential geometry, one often wishes to impose conditions
on submanifolds of a given manifold and those conditions locally look like
differential equations. Typical examples are: Lagrangian submanifolds in a
symplectic manifold, Legendrian submanifolds in a contact manifold, totally
geodesic submanifolds in a Riemannian manifold, etc. As I have already
discussed, characteristic surfaces themselves are submanifolds satisfying
suitable \textquotedblleft differential conditions\textquotedblright.
Accordingly, one speaks about \emph{PDEs imposed on submanifolds}. Jets of
sections can be generalized to \emph{jets of submanifolds}. The latter provide
a coordinate free formalism for PDEs imposed on submanifolds. In the
following, I will only need first jets of submanifolds, which can be defined
as follows.

Let $E$ be a smooth manifold. Fix a positive integer $n$ and let $\dim E=n+m$.
Consider $n$-dimensional submanifolds of $E$. Tangency at a fixed point $e\in
E$ is an equivalence relation on the set of submanifolds (through $e$). Denote
by $J_{e}^{1}(E,n)$ the set of equivalence classes. Notice that points in
$J_{e}^{1}(E,n)$ can be naturally identified with $n$-dimensional subspaces of
$T_{e}E$. Accordingly, $J_{e}^{1}(E,n)$ identifies with the Grassmannian
$\mathrm{Gr}(T_{e}E,n)$. Put
\[
J^{1}(E,n):=\coprod\limits_{e\in E}J_{e}^{1}(E,n).
\]
It identifies with the Grassmanian bundle $\mathrm{Gr}(TE,n)$. For an
$n$-dimensional submanifold $L\subset E$, denote by $[L]_{e}^{1}$ its tangency
class at $e\in L$. It is a point of $J^{1}(E,n)$ which is called the
\emph{first jet of }$L$\emph{ at }$e$. Notice that if $L_{1}$ and $L_{2}$ are
$n$-dimensional submanifolds of $E$ through the same point $e$, then there is
a (divided) chart $(x^{1},\ldots,x^{n},\boldsymbol{u})$ on $E$ which is
adapted to both, i.e., such that, in local coordinates,
\[
L_{1,2}:\boldsymbol{u}=\boldsymbol{f}_{1,2}(x^{1},\ldots,x^{n}),
\]
for some functions $\boldsymbol{f}_{1,2}=\boldsymbol{f}_{1,2}(x^{1}%
,\ldots,x^{n})$ of the $x$'s. Moreover, $L_{1}$ and $L_{2}$ have the same jet
at $e$, i.e., are tangent at $e\equiv(x_{0}^{1},\ldots,x_{0}^{n}%
,\boldsymbol{u}_{0})$, iff:
\begin{align*}
\boldsymbol{f}_{1}(x_{0}^{1},\ldots,x_{0}^{n}) &  =\boldsymbol{f}_{2}%
(x_{0}^{1},\ldots,x_{0}^{n})=\boldsymbol{u}_{0}\\
\frac{\partial\boldsymbol{f}_{1}}{\partial x^{i}}(x_{0}^{1},\ldots,x_{0}^{n})
&  =\frac{\partial\boldsymbol{f}_{2}}{\partial x^{i}}(x_{0}^{1},\ldots
,x_{0}^{n}),\quad i=1,\ldots,n.
\end{align*}
In this sense first jets of submanifolds are a coordinate free version of
first order Taylor polynomials of submanifolds. Using charts adapted to
submanifolds one can coordinatize $J^{1}(E,n)$ in an obvious way. I leave the
details to the reader. An $n$-dimensional submanifold $L$ of $E$ can be
prolonged to an $n$-dimensional submanifold $L^{(1)}$ of $J^{1}(E,n)$ by
putting
\[
L^{(1)}=\{[L]_{e}^{1}:e\in L\}.
\]
I leave to the reader to check that $L^{(1)}$ is a coordinate free version of
\textquotedblleft partial derivatives of $L$ up to the order $1$%
\textquotedblright. First jets of submanifolds are equipped with a Cartan
distribution playing the same role as in the previous subsection. A
\emph{system of first order PDEs imposed on }$n$\emph{-dimensional
submanifolds of }$E$ is a submanifold $\mathcal{E}\subset J^{1}(E,n)$. A
\emph{solutions of }$\mathcal{E}$\emph{ }is an $n$-dimensional submanifolds
$L$ of $E$ such that $L^{(1)}\subset\mathcal{E}$.

Finally, notice that if $E$ has the structure of a bundle $\pi
:E\longrightarrow M$ over an $n$-dimensional manifold $M$, then $J^{1}\pi$ is
an open and dense submanifold in $J^{1}(E,n)$.

\subsection{Singular Solutions}

\subsubsection{Multi-valued sections}

Solutions with singularities (e.g., shock waves) may have physical meaning.
For instance, in field theory, charges are often interpreted as singularities
of the fields. Therefore, it is interesting from both a mathematical and
physical point of view, to study how do singularities of solutions propagate.
We already mentioned some facts about the propagation of singularities of
solutions in the first section. Here I show that certain kinds of
singularities can be effectively treated in geometric terms within the jet
space approach to PDEs (\cite{v73,k76,k77,l80,l81,l85,v87,l90}).

Let $\pi:E\longrightarrow M$ be a fiber bundle as above, and $L$ an
$n$-dimensional, locally maximal integral submanifold of the Cartan
distribution on $J^{k}$. It is easy to see that $L$ is \emph{almost
everywhere} horizontal with respect to $\pi_{k,k-1}$ \cite{b...99}.
Consequently, $L$ is \emph{almost everywhere}, and locally, the image of a
holonomic section of $\pi_{k}$. However, $L$ doesn't need to be the image of a
holonomic section \emph{everywhere}. In particular, $L$ may project under
$\pi_{k,k-1}$, and, therefore, under $\pi_{k,0}$, to a submanifold with
singularities. Denote by $\operatorname{sing}L$ the (nowhere dense) subset of
$L$ where the singularity occur, i.e.,
\[
\operatorname{sing}L:=\{\theta\in L:d_{\theta}(\pi_{k,k-1}|_{L})\text{ is
\emph{not} injective}\}.
\]
The subset $\operatorname{sing}L\subset L$ will be referred to as
\emph{singularity locus of }$L$. A tangent space $T_{\theta}L$ to $L$ at a
point $\theta\in\operatorname{sing}L$ is called a \emph{singular }%
$R$\emph{-plane}. Singular $R$-planes may be characterized in terms of the
\emph{metaplectic structure on }$\mathcal{C}$. Namely, the correspondence
\[
\Omega:\mathcal{C\times C}\ni(X,Y)\longmapsto\Omega(X,Y):=[X,Y]+\mathcal{C}\in
TJ^{k}/\mathcal{C}%
\]
is a well-defined bilinear map (the metaplectic structure). A subspace $V$ of
$\mathcal{C}_{\theta}$, $\theta\in J^{k}$, is \emph{isotropic} iff $\Omega
(\xi,\eta)=0$ for all $\xi,\eta\in V$. (Singular) $R$-planes are
$n$-dimensional isotropic subspaces $V$. If $(d_{\theta}\pi_{k,k-1})|_{V}$ is
not injective, then $V$ is singular. A typical example of singular section is
the following: let $n=m=k=1$. In $J^{1}$ consider the smooth submanifold
\[
L:\left\{
\begin{array}
[c]{l}%
u^{2}-x^{3}=0\\
u_{x}^{2}-\dfrac{9}{4}x=0
\end{array}
\right.  .
\]
It is easy to see that $L$ is a locally maximal integral submanifold of the
Cartan distribution. However, the singularity locus of $L$ is
\[
\operatorname{sing}L=(x=0,u=0,u_{x}=0)\neq\varnothing.
\]
The projection of $L$ to $J^{0}$ is the subset
\[
L_{0}:u^{2}-x^{3}=0
\]
which has a singularity in the origin and may be interpreted as the image of
the multi-valued section.
\[
\sigma:u=\pm x^{3/2}%
\]
The other way round, the multivalued section $\sigma$ possesses a singularity
in the origin, but the singularity is \emph{resolved after the first jet
prolongation}. More generally, the multi-valued section $\sigma:u=\pm
x^{k+1/2}$ possesses a singularity in the origin which is resolved after the
$k$-th jet prolongation.

Now, let $n,m,k$ be arbitary. The above considerations suggest the following
definition: \emph{a multivalued (or singular) section of }$\pi$ \emph{is an
}$n$\emph{-dimensional, locally maximal integral submanifold of }$\mathcal{C}$.

\subsubsection{Multi-valued solutions of PDEs}

Similarly, let $\mathcal{E}\subset J^{k}$ be a PDE. Then \emph{a multivalued
solution of }$\mathcal{E}$ \emph{is an }$n$\emph{-dimensional, locally maximal
integral submanifold of }$\mathcal{C}(\mathcal{E})$. Notice that singularities
of a multivalued section $L$ are \emph{not }singularities of the submanifold
$L$ (which is always assumed to be smooth). Rather they are singularities of
the smooth map of manifolds $\pi_{k,k-1}:L\longrightarrow J^{k-1}$.
Singularities of smooth maps are usually classified along the
\emph{Thom-Boardman theory }\cite{gg73}. Here, I will only consider the
simplest one among Thom-Boardman singularities. Namely, I assume that
$d\pi_{k,k-1}|_{L}$ has constant rank $r$ along $\operatorname{sing}L$, and
that $\operatorname{sing}L\subset L$ is a smooth submanifold transversal to
$\ker(d\pi_{k,k-1}|_{L})$. In particular,
\[
\dim\operatorname{sing}L=r.
\]
It follows that the projections of $\operatorname{sing}L$ on lower order jets
are also smooth submanifolds. Let
\[
\operatorname{type}\operatorname{sing}L=n-r.
\]
For $\operatorname{type}\operatorname{sing}L=1$, $\dim\operatorname{sing}%
L=n-1$ and one speaks about \emph{fold-type singularities}. Tangent spaces to
multivaled sections with fold-type singularities at points of their singular
locus are called\emph{ type }$1$\emph{ singular }$R$\emph{-planes}. Fold-type
singularities of solutions are intimately related with characteristics. In the
following, I will only consider fold-type singularities (see, for instance,
\cite{b09} and references therein).

\subsection{Fold-type Singularities}

\subsubsection{Shapes of fold-type singularities}

Let $\mathcal{E}\subset J^{k}$ be a PDE, and $L\subset J^{k}$ a multivalued
section with a fold-type singularity along the singular locus
$\operatorname{sing}L$. Moreover, let $\theta\in\operatorname{sing}L$, and let
$S:=T_{\theta}L$ be the type $1$ singular $R$-plane tangent to $L$ at $\theta
$. If $S$ is tangent to $\mathcal{E}$ then, in a sense, $L$ is a multivalued
solution of $\mathcal{E}$ up to the order $1$. Notice that $S$, being tangent
to $\mathcal{E}$, cannot be arbitrary. Put $\underline{\theta}:=\pi
_{k,k-1}(\theta)$. Clearly,
\[
\underline{S}:=(d_{\theta}\pi_{k,k-1})(S)=(d_{\theta}\pi_{k,k-1})(T_{\theta
}\operatorname{sing}L)
\]
is an $n-1$ dimensional subspace of $T_{\underline{\theta}}J^{k-1}$. As such
it can be understood as a point in $J^{1}(J^{k-1},n-1)$. Define
\begin{align*}
\Sigma_{1}\mathcal{E} &  :=\{\underline{S}:S\text{ is a type }1\text{ singular
}R\text{-plane tangent to }\mathcal{E}\text{ at }\theta\in\mathcal{E}\}\\
&  \subset J^{1}(J^{k-1},n-1).
\end{align*}
It can be interpreted as a first order PDE for $n-1$ dimensional submanifolds
of $J^{k-1}$. Notice that if $L$ is a multivalued solution of $\mathcal{E}$
with fold-type singularity, then $\pi_{k,k-1}(\operatorname{sing}L)$ is a
solution of $\Sigma_{1}\mathcal{E}$. In this sense, $\Sigma_{1}\mathcal{E}%
$\emph{ describes the \textquotedblleft shape\textquotedblright\ of fold-type
singularities of solutions of }$\mathcal{E}$.

More precisely, let $(\boldsymbol{x},t)$ be (divided) coordinates on $M$, and,
as in Section \ref{Sec1}, denote by $\boldsymbol{u}_{\ell,J}$, $\ell+|J|{}<k$,
coordinates on $J^{k-1}$ corresponding to partial derivatives $\dfrac
{\partial^{|J|+\ell}}{\partial\boldsymbol{x}^{J}\partial t^{\ell}}$. Thus, the
$\boldsymbol{u}_{\ell,J}$ may be interpreted as derivatives along the initial
surface $t=0$ of the initial data $\partial^{\ell}\boldsymbol{u}/\partial
t^{\ell}$, $\ell<k$. Search for solutions $N$ of $\Sigma_{1}\mathcal{E}$ in
the form
\[
N:\left\{
\begin{array}
[c]{l}%
t=\tau(\boldsymbol{x})\\
\boldsymbol{u}_{\ell,J}=\tau_{\ell,J}(\boldsymbol{x}),\quad|I|{}<k
\end{array}
\right.  .
\]
Then, in general, $\Sigma_{1}\mathcal{E}$ constraints both $\tau$ and
$\tau_{\ell,J}$, looks locally like
\[
\Sigma_{1}\mathcal{E}:\boldsymbol{f}\left(  \boldsymbol{x},\ldots
,\dfrac{\partial\tau}{\partial\boldsymbol{x}^{i}},\ldots,\dfrac{\partial
\tau_{\ell,J}}{\partial\boldsymbol{x}^{i}},\ldots\right)  =0,
\]
and can therefore be interpreted as a PDE for the Cauchy data, i.e., the datum
of 1) a Cauchy surface $\Sigma:t=\tau(\boldsymbol{x})$ together with 2)
initial data $\partial^{\ell}\boldsymbol{u}/\partial t^{\ell}=\tau
_{\ell,\varnothing}(\boldsymbol{x})$, $\ell<k$, on it. When $\mathcal{E}$ is a
determined system of quasi-linear equations, then a Cauchy surface can only be
part of a solution of $\Sigma_{1}\mathcal{E}$ if it is a characteristic
surface of $\mathcal{E}$ (see below). This result relates the theory of
multivalued solutions and the theory of characteristic surfaces.

\subsubsection{The Symbol of a Differential Equations\label{SecSymb}}

Now, I want to relate fold-type singularities of solutions with
characteristics of a PDE. It will be useful to have an intrinsic definition of
characteristic covectors for a generic system of (generically fully nonlinear)
PDEs. I will present the new definition in the next section. It will
generalize (and, to some extent, clarify) the analytic definition given for
determined, quasi-linear systems. Here I provide some geometric preliminaries.

The bundle $\pi_{k,k-1}:J^{k}\longrightarrow J^{k-1}$, $k>0$, is actually an
\emph{affine bundle modelled over the vector bundle }$S^{k}T^{\ast}%
M\otimes_{J^{k-1}}VE\longrightarrow J^{k-1}$ whose fiber at $\underline
{\theta}\in J^{k-1}$ is the vector space $S^{k}T_{x}^{\ast}M\otimes V_{e}E$,
$e=\pi_{k-1,0}(\underline{\theta})$, $x=\pi_{k-1}(\underline{\theta})=\pi(e)$
(here $V_{e}E$ is the $\pi$-vertical tangent bundle to $E$). In local
coordinates, the affine structure in $\pi_{k,k-1}^{-1}(\underline{\theta})$
looks as follows. Let $\theta\in J^{k}$ and $\pi_{k,k-1}(\theta)=\underline
{\theta}$, and let $\theta$ have jet coordinates
\[
(x^{1},\ldots,x^{n},\ldots,\boldsymbol{u}_{I},\ldots)\equiv\theta,\quad
|I|{}\leq k.
\]
Take $v\in S^{k}T_{x}^{\ast}M\otimes V_{e}E$, and let
\[
v=\boldsymbol{v}_{i_{1}\cdots i_{k}}dx^{i_{1}}\cdots dx^{i_{k}}\otimes
\dfrac{\partial}{\partial\boldsymbol{u}}.
\]
One can use the $\boldsymbol{v}_{i_{1}\cdots i_{k}}$'s as coordinates in
$S^{k}T_{x}^{\ast}M\otimes V_{e}E$. Then $\theta+v\in\pi_{k,k-1}%
^{-1}(\underline{\theta})$ have jet coordinates
\begin{equation}
(x^{1},\ldots,x^{n},\ldots,\boldsymbol{u}_{J},\ldots,\boldsymbol{u}%
_{i_{1}\cdots i_{k}}+\tfrac{k!}{(i_{1}\cdots i_{k})!}\boldsymbol{v}%
_{i_{1}\cdots i_{k}},\ldots)\equiv\theta+v,\quad|J|{}<k,\label{4}%
\end{equation}
(see Subsection \ref{SecNL} for the meaning of the combinatorial coefficient
$(i_{1}\cdots i_{k})!/k!$). As a consequence, the vertical bundle
$VJ^{k}\longrightarrow J^{k}$ of $\pi_{k,k-1}$ is isomorphic to the vector
bundle $J^{k}\times_{J^{k-1}}\Sigma_{k-1}\longrightarrow J^{k}$. In local
coordinates, the isomorphism looks as
\[
\dfrac{\partial}{\partial\boldsymbol{u}_{i_{1}\cdots i_{k}}}\longmapsto
\frac{(i_{1}\cdots i_{k})!}{k!}\,dx^{i_{1}}\cdots dx^{i_{k}}\otimes
\dfrac{\partial}{\partial\boldsymbol{u}}.
\]
Let $\mathcal{E}\subset J^{k}$ be a PDE locally given by (\ref{5}). According
to (\ref{4}), $\mathcal{E}$\emph{ has a quasi-linear local description iff it
is an affine subbundle of }$J^{k}\longrightarrow J^{k-1}$. This remark
provides an intrinsic definition of quasi-linear equations. Now, let
$\mathcal{E}$ be generic, $\theta\in\mathcal{E}$, $\underline{\theta}%
:=\pi_{k,k-1}(\theta),$ $e=\pi_{k,0}(\theta)$, and $x=\pi_{k}(\theta)$. Put
$g_{\theta}:=V_{\theta}J^{k}\cap T_{\theta}\mathcal{E}$. In view of the affine
structure in the fibers of $\pi_{k,k-1}$, $g_{\theta}$ can be understood as a
subspace of $S^{k}T_{x}^{\ast}M\otimes V_{e}E$. It is easy to see that
$g_{\theta}$ consists of $v\in S^{k}T_{x}^{\ast}M\otimes V_{e}E$, such that
\[
g_{\theta}:\sum_{i_{1}\leq\cdots\leq i_{k}}\,\dfrac{\partial\boldsymbol{F}%
}{\partial\boldsymbol{u}_{i_{1}\cdots i_{k}}}\cdot\boldsymbol{v}_{i_{1}\cdots
i_{k}}=0
\]
In particular, if $\mathcal{E}$ is determined, and quasi-linear then
\[
g_{\theta}:\boldsymbol{A}^{i_{1}\cdots i_{k}}\cdot\boldsymbol{v}_{i_{1}\cdots
i_{k}}=0,
\]
where $\boldsymbol{A}=(\boldsymbol{A}^{i_{1}\cdots i_{k}})$ is the symbol of
$\boldsymbol{F}$. For this reason $g_{\theta}$ is called the \emph{symbol of
}$\mathcal{E}$\emph{ at }$\theta$. In the case when $\mathcal{E}$ is a linear
equation $\boldsymbol{Du}=0$, $\boldsymbol{D}$ being a linear differential
operator, the symbol $g_{\theta}$ does only depend on $x=\pi_{k}(\theta)$. If,
moreover, $M$ is an Euclidean space, the symbol can be understood as a
homogeneous fiber-wise polynomial function $\sigma(\boldsymbol{D})$ on
$T^{\ast}M$. In this case, it can be defined analytically as
\[
\sigma(\boldsymbol{D})(\boldsymbol{p})=0,\quad\sigma(\boldsymbol{D}%
)(\boldsymbol{p}):=e^{-\boldsymbol{p}\cdot\boldsymbol{x}}\boldsymbol{D}%
(e^{\boldsymbol{p}\cdot\boldsymbol{x}}),\quad\boldsymbol{x}\in M,\quad
\boldsymbol{p}\in T_{\boldsymbol{x}}M,
\]
and plays an important role in quantization (see, e.g., \cite{ems04}).

\subsubsection{Characteristic Covectors of a PDE}

Let $\mathcal{E}$, $\theta$, $e$, $x$ be as in the above subsection. A non
zero covector $\boldsymbol{p}\in T_{x}^{\ast}M$ is a \emph{characteristic
covector} for $\mathcal{E}$ at $\theta$ if there exists a non zero $\xi\in
V_{e}E$ such that
\[
\boldsymbol{p}\cdot\cdots\cdot\boldsymbol{p}\otimes\xi\in g_{\theta}.
\]
If $\boldsymbol{p}=p_{i}dx^{i}$ and $\xi=\boldsymbol{\xi\,}\partial
/\partial\boldsymbol{u}$, this means that the system of linear equations
\[
\sum_{i_{1}\leq\cdots\leq i_{k}}\,p_{i_{1}}\cdots p_{i_{k}}\left.
\dfrac{\partial\boldsymbol{F}}{\partial\boldsymbol{u}_{i_{1}\cdots i_{k}}%
}\right\vert _{\theta}\cdot\boldsymbol{\xi}%
\]
in the $m$ unknowns $\boldsymbol{\xi}$ has non-trivial solutions. In other
words
\begin{equation}
\operatorname{rank}\sum_{i_{1}\leq\cdots\leq i_{k}}\,p_{i_{1}}\cdots p_{i_{k}%
}\left.  \dfrac{\partial\boldsymbol{F}}{\partial\boldsymbol{u}_{i_{1}\cdots
i_{k}}}\right\vert _{\theta}<m.\label{7}%
\end{equation}
Notice that for underdetermined systems, i.e., $\operatorname{codim}%
\mathcal{E}\leq m$, every covector is characteristic, and, therefore, only the
determined and overdetermined cases ($\operatorname{codim}\mathcal{E}\geq m$)
need to be considered. For a determined system, condition (\ref{7}) coincides
with condition
\[
\det\sum_{i_{1}\leq\cdots\leq i_{k}}\,p_{i_{1}}\cdots p_{i_{k}}\left.
\dfrac{\partial\boldsymbol{F}}{\partial\boldsymbol{u}_{i_{1}\cdots i_{k}}%
}\right\vert _{\theta}=0
\]
which I already considered in Section \ref{Sec1}. For quasi-linear, determined
systems, $\boldsymbol{p}=p_{i}dx^{i}$ is characteristic iff
\[
\det\boldsymbol{A}(\boldsymbol{p})=0.
\]
One concludes that the notion introduced here is a coordinate free version of
the one introduced in Section \ref{Sec1}. Notice that the characteristic
condition (\ref{7}) depends on the point $\theta$ in $\mathcal{E}$ and, in
general, changes from point to point.

\subsubsection{Characteristic Covectors and Fold-type Singularities}

In Subsection \ref{SecPhys}, I presented an (informal) argument showing that
characteristics are actually related to singularities of solutions: namely,
singularities of solutions of quasi-linear, determined systems of PDEs occur
along characteristic surfaces. Here I present a rigorous, intrinsic argument
which applies to generic nonlinear PDEs. First I need few remarks. Let
$\theta=[\sigma]_{x}^{k}\in\mathcal{E}$. The set of characteristic covectors
at $\theta$ is the \emph{characteristic variety} at $\theta$. In this way, one
gets a (possibly singular) fiber bundle over $\mathcal{E}$ whose fibers are,
by definition, characteristic varieties. Call it the \emph{characteristic
bundle of }$\mathcal{E}$. Now, consider the $R$-plane $R_{\theta}\subset
T_{\underline{\theta}}J^{k-1}$ corresponding to $\theta$. The projection
$d_{\underline{\theta}}\pi_{k-1}$ identifies $R_{\theta}$ with $T_{x}M$. Thus,
if $\boldsymbol{p}$ is a non-zero covector in $T_{x}^{\ast}M$, one can
understand its kernel as an $(n-1)$-dimensional subspace in $R_{\theta}$. I
denote it by $\ker_{\theta}\boldsymbol{p}$. Specifically, $\ker_{\theta
}\boldsymbol{p}:=(d_{x}j^{k-1}\sigma)(\ker\boldsymbol{p})$. It holds the
following proposition: \emph{if }$\mathcal{E}$\emph{ is a formally integrable
PDE then the equation of fold-type singularities is \textquotedblleft
dual\textquotedblright\ to the characteristic bundle in the following sense: }%
\begin{equation}
\Sigma_{1}\mathcal{E}=\{\ker_{\theta}\boldsymbol{p}:\boldsymbol{p}\text{ is a
characteristic covector of }\mathcal{E}\text{ at }\theta\in\mathcal{E}%
\}.\label{8}%
\end{equation}
(for this part of the statement see, for instance, \cite[Theorem 5.2]{b09},
Theorem 5.2. Formal integrability roughly means that: if a Taylor polynomial
of the order $k+\ell$ is a solution of $\mathcal{E}$ up to the order $\ell$,
then it can be \textquotedblleft completed\textquotedblright\ to a
\emph{formal solution}, i.e., a solution in the form of a (possibly
non-converging) Taylor series. \emph{If, moreover, }$\mathcal{E}$\emph{ is
determined, quasi-linear and }$A(\boldsymbol{p})$\emph{ is generically
invertible (i.e., }$\mathcal{E}$\emph{ can be generically recast in normal
form), then }$\Sigma_{1}E$\emph{ is locally equivalent to (\ref{CEdiv}) (which
constraints the shape of a Cauchy surface }$\Sigma$\emph{) + (\ref{Sigma})
(which constraints the initial data on }$\Sigma$\emph{) }(see \cite{v87}).
Concluding, the equation for fold-type singularities is an equation for
\emph{Cauchy data} (including both the Cauchy surface and the initial data on
it), telling us that 1) fold-type singularities may only occur along a
characteristic surface and that 2) initial data on a characteristic surface
may not be assigned arbitrarily.

\subsubsection{An Example: the 2D Klein-Gordon Equation}

The equation for characteristic surfaces of a quasi-linear system does not
contain a full information on the original equation. In fact, it does only
depend on its symbol. However, the equation for singularities may contain a
full information on the original equation. I briefly illustrate this
phenomenon with a simple example. For details, see \cite{v73,v87} (see also
\cite{l...94}, where more examples from Mathematical Physics can be found).

Consider the Klein-Gordon equation on the $2$-dimensional Minkowski
space-time:
\[
\mathcal{E}_{KG}:u_{tt}-u_{\boldsymbol{xx}}+m^{2}u=0.
\]
I already showed that characteristic surfaces of the Klein-Gordon equation are
null hypersurfaces. An hypersurface $\Sigma:z(t,x)=0$ is a characteristic
surface iff
\[
(z_{t}^{2}-z_{x}^{2})|_{z=0}=0
\]
This shows that $z_{t}\neq0$ so that $\Sigma$ is actually of the form
$\Sigma:t=\tau(x)$, with
\[
\tau_{x}^{2}=1
\]
which is the $1$-dimensional eikonal equation. The fold-type singularity
equation $\Sigma_{1}\mathcal{E}_{KG}\subset J^{1}(J^{1},1)$ can be easily
computed, using, for instance, (\ref{8}). Coordinatize $J^{1}$ by
$x,t,u,u_{x},u_{t}$. Since characteristic surfaces are of the form
$\Sigma:t=\tau(x)$ we can interpret $x$ as independent variable, and
coordinatize $J^{1}(J^{1},1)$ by $x,t,u,u_{x},u_{t},t^{\prime},u^{\prime
},u_{x}^{\prime},u_{t}^{\prime}$, where $f^{\prime}$ means $df/dx$. Then
\[
\Sigma_{1}\mathcal{E}_{KG}:\left\{
\begin{array}
[c]{l}%
(t^{\prime})^{2}=1\\
u^{\prime}=u_{x}+t^{\prime}u_{t}\\
u_{x}^{\prime}=m^{2}u+t^{\prime}u_{t}^{\prime}%
\end{array}
\right.  .
\]
Eliminating $u_{x}$, one gets the following second order system for the Cauchy
data $t,u,u_{t}$:
\[
\left\{
\begin{array}
[c]{l}%
(t^{\prime})^{2}=1\\
u^{\prime\prime}-t^{\prime\prime}u_{t}-2t^{\prime}u_{t}^{\prime}=m^{2}u
\end{array}
\right.  .
\]
Notice that $\Sigma_{1}\mathcal{E}_{KG}$ contains the mass parameter $m$.
Actually, it can be proved, by purely geometric methods, that $\Sigma
_{1}\mathcal{E}_{KG}$ contains a full information about $\mathcal{E}_{KG}$.
More generally, understanding \emph{when a PDE can be reconstructed from its
singularity equation} is an interesting open problem (in some sense, as I
already outlined in the introduction, analogous to quantization) that has been
first addressed by Vinogradov in simple situations \cite{v73,v87}.

\subsubsection{Bicharacteristics of Determined Systems of PDEs}

Let $\mathcal{E}\subset J^{k}$ be a determined system of PDEs, and $L$ a
multi-valued solution with a fold-type singularity along $\operatorname{sing}%
L$. If one interprets $L$ as a wave propagating in the space-time, then it is
natural to intepret $\pi_{k}(\operatorname{sing}L)\subset M$ as its
\emph{wave-front}. Recall that the wave-front is a characteristic surface. It
can be shown that $\operatorname{sing}L$ is equipped with a canonical field of
directions, i.e., a $1$-dimensional distribution. Accordingly, the wave-front
of $L$ is foliated by $1$-dimensional submanifolds \cite{k77,klv86}. In the
case when $\mathcal{E}$ is a linear system, this is a classical result, the
$1$-dimensional leaves of $L$ are called \emph{bicharacteristics}, and one
usually says that \emph{wave-fronts propagate along bicharacteristics}. I will
not discuss this result in full generality, which would require too much
space. Instead, I will consider, in the next section, the case when the symbol
$g_{\theta}$ at $\theta\in\mathcal{E}$ does only depend on $x=\pi_{k}(\theta
)$. In this case the wave-front is a solution of a genuine first order PDE in
one dependent variable, which locally looks like (\ref{CEdiv}), where the
$B$'s does only depend on the $\boldsymbol{x}$'s and $\tau$ . We are thus led
to consider the class of first order PDEs in one dependent variable. This will
be the main topic of the last section. Notice once again that, from a physical
point of view, the passage from $\pi_{k}(\operatorname{sing}L)$ to its
bicharacteristics can be interpreted as the passage from a wave optics (the
dynamics of wave-fronts) to a geometric optics (the dynamics of rays).

\section{Bicharacteristics and the Hamilton-Jacobi Theory\label{SecBic}}

\subsection{Contact Geometry of First Jets of Functions}

\subsubsection{Jets of functions}

In this section, I focus on first order PDEs in one dependent variable. The
equations for characteristic surfaces of determined, quasi-linear systems
whose symbol does only depend on independent variables (for instance, linear
systems) are precisely of this kind. For the sake of simplicity, I will
suppose that the equation under consideration is imposed on a real function
$f$ on a manifold $M$. This is always true locally. More generally, one could
consider equations imposed on sections of a bundle with one dimensional fibers
or on $1$-codimensional submanifolds of a given manifold. Similar results as
the one presented in this lecture hold for these (unparameterized) cases.

Notice that a real function $f$ on a manifold $M$ can be understood as a
section of the trivial bundle $\pi_{M}:M\times\mathbb{R}\longrightarrow M$. In
the following, the first jet space of $\pi_{M}$ will be denoted by $J^{1}(M)$.
A \emph{first order PDE in one dependent variable} is then a hypersurface
$\mathcal{E}$ in $J^{1}(M)$. The most remarkable property of $J^{1}(M)$ is
that \emph{it is equipped with a natural contact structure}. Recall that a
\emph{contact structure} on a $2n+1$ dimensional manifold $N$ is an hyperplane
distribution with non-degenerate, associated metaplectic structure. A contact
structure can be presented as the kernel distribution $\ker\alpha$ of a
\emph{contact form} $1$-form $\alpha$ such that $d\alpha$ is non degenerate on
$\ker\alpha$. \emph{The contact structure in }$J^{1}(M)$\emph{ is given by the
Cartan distribution}. A contact $1$-form on $J^{1}(M)$ can be defined as
follows. First of all, notice that there is a canonical isomorphism (of
bundles over $M$)
\[
J^{1}(M)\simeq T^{\ast}M\times\mathbb{R}%
\]
given by $[f]_{x}^{1}\longmapsto(d_{x}f,f(x))$. Denote by $u:M\times
\mathbb{R}\longrightarrow\mathbb{R}$ the canonical function on $M\times
\mathbb{R}$, i.e., the projection onto the second factor, and by $\theta$ the
tautological $1$-form on $T^{\ast}M$. Abusing the notation, I denote by the
same symbols $u$ and $\theta$, the pull-backs on $J^{1}(M)$. The $1$-form
\[
\alpha:=du-\theta\in\Lambda^{1}(J^{1}(M))
\]
is a contact form. Indeed, it is easy to see that its coordinate description
in jet coordinates is
\[
\alpha=du-u_{i}dx^{i}.
\]
Moreover, $\ker\alpha$ is precisely the $2n$-dimensional Cartan distribution,
and $d\alpha$ is non-degenerate over it. The contact geometry of $J^{1}(M)$ is
intimately related to the symplectic geometry of $T^{\ast}M$.

\subsubsection{Jacobi Algebra of a Contact Manifold}

Recall that functions on a symplectic manifold form a Poisson algebra equipped
with a morphism of Lie algebras into infinitesimal syplectomorphisms.
Similarly, functions on a contact manifold form a Jacobi algebra equipped with
a morphism of Lie algebras into \emph{infinitesimal contactomorphisms}. Let us
illustrate this in the simple case of the contact manifold $J^{1}(M)$. In this
case, a \emph{contactomorphism} is nothing but a diffeomeorphism
$J^{1}(M)\longrightarrow J^{1}(M)$ preserving the Cartan distribution.
Similarly, an \emph{infinitesimal contactomorphism} is a vector field $X$ over
$J^{1}(M)$ whose flow preserves the Cartan distribution. In other words,
\[
L_{X}\alpha=\lambda\alpha\quad\text{for some function }\lambda\in C^{\infty
}(J^{1}(M)).
\]
A smooth function $f$ on $J^{1}(M)$ determines an infinitesimal
contactomorphism and vice-versa as follows. Let $\partial/\partial u$ be the
vector field on $J^{1}(M)$ determined by the canonical coordinate vector field
on $\mathbb{R}$ and the identification $J^{1}(M)\simeq T^{\ast}M\times
\mathbb{R}$. The vector field $\partial/\partial u$ is transversal to the
Cartan distribution, so that $TJ^{1}(M)=\langle\partial/\partial
u\rangle\otimes\mathcal{C}$. Consider the $1$-form
\[
\delta f:=df-\frac{\partial f}{\partial u}\alpha.
\]
It is easy to see that $\delta f\in\operatorname{Ann}(\partial/\partial u)$ so
that there exists a unique vector field $Y_{f}$ in the Cartan distribution
such that
\[
i_{Y_{f}}d\alpha=\delta f.
\]
The vector field
\[
X_{f}:=Y_{f}-f\frac{\partial}{\partial u}%
\]
is an infinitesimal contactomorphism and every infinitesimal contactomorphism
$X$ is of the form $X=X_{f}$, with $f=-\alpha(X)$. Notice that $\partial
/\partial u=-X_{1}$. Finally, for any two smooth functions $f,g$ on $J^{1}(M)$
one has
\begin{align*}
X_{fg} &  =fX_{g}+gX_{g}-fgX_{1}\\
\lbrack X_{f},X_{g}] &  =X_{\{f,g\}}%
\end{align*}
with
\[
\{f,g\}:=X_{f}(g)-X_{1}(f)g
\]
This shows that smooth functions on $J^{1}(M)$ equipped with the bracket
$\{-,-\}$ form a Jacobi algebra isomorphic (as a Lie algebra) to the Lie
agebra of infinitesimal contactomorphisms. In local coordinates,
\[
X_{f}=\frac{\partial f}{\partial u_{i}}\frac{\partial}{\partial x^{i}}-\left(
\frac{\partial f}{\partial x^{i}}+u_{i}\frac{\partial f}{\partial u}\right)
\frac{\partial}{\partial u_{i}}+\left(  u_{i}\frac{\partial f}{\partial u_{i}%
}-f\right)  \frac{\partial}{\partial u}%
\]
and
\[
\{f,g\}=\frac{\partial f}{\partial u_{i}}\frac{\partial g}{\partial x^{i}%
}-\frac{\partial g}{\partial u_{i}}\frac{\partial f}{\partial x^{i}}%
+u_{i}\left(  \frac{\partial f}{\partial u_{i}}\frac{\partial g}{\partial
u}-\frac{\partial g}{\partial u_{i}}\frac{\partial f}{\partial u}\right)
-f\frac{\partial g}{\partial u}+g\frac{\partial f}{\partial u}.
\]

\subsection{First Order Scalar PDEs}

\subsubsection{(Bi)characteristic Foliation and the Method of Characteristics}

Now, let $\mathcal{E}\subset J^{1}(M)$ be a (codimension $1$) PDE. As a
minimal regularity condition, I assume that the Cartan distribution
$\mathcal{C}(\mathcal{E})$ on $\mathcal{E}$ is regular (i.e., constant
dimension). It then follows that
\[
\dim\mathcal{C}(\mathcal{E})=2n-1.
\]
Since $d\alpha$ is a symplectic form on $\mathcal{C}$, it must degenerate on
$\mathcal{C}(\mathcal{E})$ along a field of directions $\ell(\mathcal{E}%
)\subset\mathcal{C}(\mathcal{E})$ on $\mathcal{E}$. Integral manifolds of
$\ell(\mathcal{E})$ foliate $\mathcal{E}$ and are called \emph{characteristic
lines} of $\mathcal{E}$, or, \emph{bicharacteristics }if $\mathcal{E}$ is the
equation for characteristic surfaces of a determined, quasi-linear system
(whose symbol does only depend on independent variables). If $\mathcal{E}$ is
assigned as the zero locus of a function $F$ on $J^{1}(M)$, i.e.,
$\mathcal{E}:F=0$, then $\ell(\mathcal{E})$ is spanned by the vector field
$Y_{F}$.

The key remark here is that $\ell(E)$\emph{ is tangent to every (multivalued)
solution of }$E$. Therefore, solutions themselves are foliated by
$1$-dimensional leaves. More generally, it can be proved that \emph{solutions
of the fold-type singularity equation of a determined system of PDEs are
foliated by }$1$\emph{-dimensional leaves}. One concludes that
\emph{singularities of solutions of determined systems of PDEs propagate along
bicharacteristics}.

The existence of characteristic lines suggests a way to solve the Cauchy
problem for $\mathcal{E}$. Namely, let $\Sigma\subset M$ be an hypersurface
and $\mu$ a smooth function on it. Search for a solution $f$ of $\mathcal{E}$
such that $f|_{\Sigma}=\mu$, i.e., understand $(\Sigma,\mu)$ as Cauchy data
for $\mathcal{E}$. If $\Sigma$ is not a characteristic surface for
$\mathcal{E}$, then the goal can be achieved as follows. First, notice that
there exists a unique $(n-1)$-dimensional submanifold $N$ in $\mathcal{E}$
such that 1) $N$ projects to the graph of $\mu$ under $\pi_{1,0}$, 2) $N$ is
integral for $\mathcal{C}(\mathcal{E})$, 3) $N$ is transversal to
$\ell(\mathcal{E})$ \cite{b...99} (see also \cite{vk81} for the case of an
Hamilton-Jacobi equation). The submanifold $N$ encodes the information about
$\Sigma$, $\mu$, and derivatives of $\mu$ along $\Sigma$. The union of
characteristic lines passing through $N$ is, by construction, an
$n$-dimensional integral manifold of $\mathcal{C}(\mathcal{E})$ horizontal
with respect to fibers of $\pi_{1,0}$. As such, it is the image of the first
jet prolongation of a solution of $\mathcal{E}$ agreeing with the Cauchy data
$(\Sigma,\mu)$. Notice that, if $\mathcal{E}:F=0$, and $\Sigma$ is in the form
$\Sigma:z=0$, then the condition of not-being characteristic is
\[
\left.  \frac{\partial F}{\partial u_{i}}\frac{\partial z}{\partial x^{i}%
}\right\vert _{z=0}\neq0.
\]
In this case, solving the assigned Cauchy problem amounts to solve the
following system of ODEs:
\[
\left\{
\begin{array}
[c]{l}%
\dot{x}^{i}=\frac{\partial F}{\partial u_{i}}\\
\dot{u}_{i}=-\frac{\partial F}{\partial x^{i}}-u_{i}\frac{\partial F}{\partial
u}\\
\dot{u}=u_{i}\frac{\partial F}{\partial u_{i}}%
\end{array}
\right.  ,
\]
with initial data on $N$. This is nothing but the classical method of
characteristics to solve $1$st order scalar PDEs.

\subsubsection{An example}

Consider the following Cauchy problem in two independent variables
$x^{1},x^{2}$:
\begin{align*}
u-\frac{\partial u}{\partial x^{1}}\frac{\partial u}{\partial x^{2}} &  =0,\\
u|_{x^{2}=0} &  =(x^{1})^{2}.
\end{align*}
Then $F=u-u_{1}u_{2}$ and $\Sigma:z(x^{1},x^{2})=x^{2}=0$ is (almost
everywhere) non-characteristic since
\[
\left.  \frac{\partial F}{\partial u_{i}}\frac{\partial z}{\partial x^{i}%
}\right\vert _{z=0}=-u_{1}|_{x^{2}=0}=-2x^{1}\neq0\quad\text{almost
everywhere}.
\]
The Cauchy data determine a $1$-dimensional integral manifold $N$ for
$\mathcal{C}(\mathcal{E})$ which is parametrically given by
\begin{equation}
N:\left\{
\begin{array}
[c]{c}%
x^{1}=s\\
x^{2}=0\\
u=s^{2}\\
u_{1}=2s\\
u_{2}=s/2
\end{array}
\right.  .\label{12}%
\end{equation}
Indeed, one may check that this is the unique choice of $N$ satisfying all the
required properties. One also has
\[
Y_{F}=-u_{2}\frac{\partial}{\partial x^{1}}-u_{1}\frac{\partial}{\partial
x^{2}}-2u_{1}u_{2}\frac{\partial}{\partial u}-u_{1}\frac{\partial}{\partial
u_{1}}-u_{2}\frac{\partial}{\partial u_{2}}.
\]
So characteristic lines may be computed integrating equations
\[
\left\{
\begin{array}
[c]{c}%
\dot{x}^{1}=-u_{2}\\
\dot{x}^{2}=-u_{1}\\
\dot{u}=-2u_{1}u_{2}\\
\dot{u}_{1}=-u_{1}\\
\dot{u}_{2}=-u_{2}%
\end{array}
\right.  ,
\]
with (parametric) initial conditions given by (\ref{12}). Integrating and
eliminating the parameters (and the higher derivatives) one gets the solution
\[
u=\frac{(4x^{1}+x^{2})^{2}}{16}.
\]

\subsection{Hamilton-Jacobi Theory}

\subsubsection{Hamilton-Jacobi equations}

It may happen that $\mathcal{E}$ is the preimage of an hypersurface
$\mathcal{H}$ of $T^{\ast}M$ under the canonical projection $J^{1}%
(M)\longrightarrow T^{\ast}M$. If $\mathcal{E}$ is locally given by
$\mathcal{E}:\{F=0$, then $F$ can be chosen such that $\frac{\partial
F}{\partial u}=0$, i.e., $F=F(x^{1},\ldots,x^{n},u_{1},\ldots,u_{n})$. In
other words, $F$ is the pull-back of a function $H=H(x^{1},\ldots,x^{n}%
,p_{1},\ldots,p_{n})$ on $T^{\ast}M$. In this case, $\mathcal{E}$ is precisely
the \emph{Hamilton-Jacobi equation }associated to the Hamiltonian system
$(T^{\ast}M,H)$. Finding characteristics of $\mathcal{E}$ is then the same as
finding \emph{characteristics of }$\mathcal{H}$, i.e., the degeneracy lines of
the restriction to $\mathcal{H}$ of the canonical symplectic form
$\Omega:=-d\theta$ on $T^{\ast}M$. In their turn, characteristics of
$\mathcal{H}$ are trajectories of the Hamiltonian vector field $X_{H}$ defined
by
\[
i_{X_{H}}\Omega=dH.
\]
Thus the method of characteristics to solve a Cauchy problem for an
Hamilton-Jacobi equation consists in integrating the Hamilton equations with
suitable initial data.

\subsubsection{Hamilton-Jacobi theorem}

I conclude this section reviewing the Hamilton-Jacobi theory of a Hamiltonian
system $(T^{\ast}M,H)$. First, I specialize the geometric definition of PDE,
(multivalued) solutions, and the method of (bi)characteristics to this context.

As already noticed, if one is interested in first order PDEs, in one
independent variable, of the form $H(x^{1},\ldots,x^{n},u_{1},\ldots,u_{n}%
)=E$, $E$ being a constant, then one may understand them geometrically as
hypersurfaces $\mathcal{H}$ in $T^{\ast}M$, $M$ being a manifold (of
independent variables) coordinatized by $x^{1},\ldots,x^{n}$. Any such
hypersurface will be referred to as a \emph{Hamilton-Jacobi equation}. The
cotangent bundle $T^{\ast}M$ comes equipped with its canonical symplectic
structure $\Omega$. Locally, $\Omega=dp_{i}\wedge dx^{i}$, where the $p_{i}$'s
are cotangent coordinates conjugate to the $x^{i}$'s. The symplectic form
$\Omega$ plays here a similar role as the Cartan distribution in the general
theory of PDEs. The geometric definition of solutions of $\mathcal{H}$ is
clear: a solution is a (local) function $f$ on $M$ such that $df$ takes values
into $\mathcal{H}$. Notice that the image of $df$ is a Lagrangian submanifold
of $T^{\ast}M$ horizontal with respect to the projection $T^{\ast
}M\longrightarrow M$, and viceversa: every Lagrangian submanifold horizontal
with respect to $T^{\ast}M\longrightarrow M$ is locally the image of $df$ for
some local function $f$ on $M$. Similarly as in the general case, one could be
interested in \emph{multivalued (singular) solutions}. As in the general case,
a geometric definition is obtained relaxing the horizontality condition from
the definition of a solution. Thus, \emph{a multivalued solution of a Hamilton
Jacobi equation }$\mathcal{H}\subset T^{\ast}M$ is a Lagrangian submanifold
$L\subset T^{\ast}M$ such that $L\subset\mathcal{H}$.

\begin{example}
Let $M=\mathbb{R}$, and $\mathcal{H}:p^{2}+q^{2}=E$, $E>0$. Then $\mathcal{H}$
is a multivalued solution, corresponding to the multivalued function
implicitly defined by
\[
\left(  \frac{df}{dq}\right)  ^{2}+q^{2}=E.
\]

\end{example}

Now notice that every Hamilton-Jacobi equation $\mathcal{H}$ comes equipped
with a canonical field of directions $\ell(\mathcal{H})$: the \emph{degeneracy
distribution of }$\Omega|_{\mathcal{H}}$. Namely, for $\boldsymbol{p}%
\in\mathcal{H}$
\[
\ell(\mathcal{H})_{\boldsymbol{p}}:=\{\xi\in T_{\boldsymbol{p}}\mathcal{H}%
:i_{\xi}\Omega|_{\mathcal{H}}=0\}.
\]
One of the key points here is the following version of the
\emph{Hamilton-Jacobi Theorem}:\emph{ }$\ell(\mathcal{H})$\emph{ is tangent to
every multivalued solution of }$\mathcal{H}$. In particular, $\ell
(\mathcal{H})$ restricts to a field of (bi)characteristic directions on any
(multivalued) solution. When $\mathcal{H}:H=E$ for some function $H$ on
$T^{\ast}M$, then the Hamiltonian vector field $X_{H}$ of $H$ is tangent to
$\mathcal{H}$ and generates the field of directions $\ell(\mathcal{H})$.
Therefore, to find (multivalued) solutions of the Hamilton-Jacobi equation
$\mathcal{H}$, it is enough to start from an $(n-1)$-dimensional submanifold
$N\subset$ $\mathcal{H}$ such that 1) $N$ is isotropic, 2) $N$ is transversal
to $X_{H}$, and then move $N$ along the flow of $X_{H}$. The $n$-dimensional
submanifold swept in this way is, by construction, Lagrangian, and, therefore,
is a multivalued solution. This procedure specializes the method of
characteristics to Hamilton-Jacobi equations. Thus, \emph{one can get
solutions of the Hamilton-Jacobi equation }$H:H=E$\emph{, from solutions of
the Hamilton equations} (for the flow of $X_{H}$). The converse is also true
to some extent: any multivalued solution $L$ of the Hamilton-Jacobi equation
is invariant under the flow of $X_{H}$: restricting $X_{H}$ to $L$ reduces by
$n$ the number of degrees of freedom in the Hamilton equations, and
simplifyies the integration problem.

The geometric version of the Hamilton-Jacobi theorem recalled above (which is
actually a simple remark in geometric terms) has been generalized by the
Geometric Mechanics community to many different contexts: Lagrangian mechanics
\cite{c...06}, non-holonomic systems \cite{idd08,dmd08b,ob09,c...10,los12},
almost Poisson manifolds \cite{ddv12}, mechanics on Lie algebroids
\cite{ls12}, field theory \cite{dmd08,d...10,dv13} and higher derivative field
theory \cite{v10,v11,v12}. Another key aspect of the Hamilton-Jacobi theory is
the role played by \emph{complete integrals of the Hamilton-Jacobi problem.}

\subsubsection{Complete integrals of the Hamilton-Jacobi problem}

Consider the family of Hamilton-Jacobi equations $H=E$. It is sometimes
collectively referred to as the \emph{Hamilton-Jacobi problem}. A (local)
\emph{complete integral of the Hamilton-Jacobi problem }is then a (local)
Lagrangian foliation $\mathcal{F}$ of $T^{\ast}M$ whose leaves are multivalued
solutions. Notice that having a complete integral amounts to having an
$n$\emph{-parameter family of solutions depending on the parameters in an
essential way}. Now, suppose that the space of leaves of $\mathcal{F}$ is a
smooth manifold $Q$ and that the canonical map $\varphi:T^{\ast}%
M\longrightarrow Q$ is a submersion. Let $q^{1},\ldots,q^{n}$ be coordinates
in $Q$. Thus, they are precisely the $n$-parameters parameterizing athe
complete integral. Interpret the $q^{i}$'s as (local) functions on $T^{\ast}%
M$. Since fibers of $\varphi$ are contained into the level surfaces of $H$,
one can always choose one of the $q^{i}$'s to be $H$ itself. Applying the
Hamilton-Jacobi theorem to the Hamilton-Jacobi equations $q^{i}=c^{i}$, with
$c^{i}$'s constant, one sees that the $q^{i}$'s are actually $n$ independent
and Poisson-commuting functions on $T^{\ast}M$ (see, for instance,
\cite{vk81}). It follows that $(T^{\ast}M,\Omega,H)$ is a (locally) integrable
Hamiltonian system! This result clarifies the use of the Hamilton-Jacobi
problem to integrate Hamilton equations.

On another hand, let $\mathcal{F}$ be a complete integral of the
Hamilton-Jacobi problem such that: 1) the space of leaves of $\mathcal{F}$ is
a smooth manifold $Q$ with local coordinates $q^{1},\ldots,q^{n}$ and the
canonical map $\varphi:T^{\ast}M\longrightarrow Q$ is a submersion, 2) the
leaves of $\mathcal{F}$ are all graphs of (closed) $1$-forms, i.e., sections
of $T^{\ast}M\longrightarrow M$. It follows thet there is a diffeomorphism
$\Phi:M\times Q\simeq T^{\ast}M$, locally given by:
\[
\Phi^{\ast}(x^{i})=x^{i},\quad\Phi^{\ast}(p_{i})=\frac{\partial W}{\partial
x^{i}},
\]
for some local function $W=W(x^{1},\ldots,x^{n},q^{1},\ldots,q^{n})$. Clearly,
$M\times Q$ inherits a symplectic structure $\Omega_{\mathcal{F}}:=\Phi_{\ast
}(\Omega)$ and an Hamiltonian function $H_{\mathcal{F}}:=\Phi_{\ast}(\Omega)$.
It is easy to see that the local functions on $M\times Q$ defined by%
\[
P_{i}:=\frac{\partial W}{\partial q^{i}}%
\]
are conjugate to the $q^{i}$'s, i.e.,
\[
\Omega_{\mathcal{F}}=dP_{i}\wedge dq^{i},
\]
moreover, the Hamiltonian system $(M\times Q,\Omega_{\mathcal{F}%
},H_{\mathcal{F}})$ is canonically isomorphic to $(T^{\ast}M,\Omega,H)$, but
$H_{\mathcal{F}}$ does not depend on the $P_{i}$'s. Therefore, in the
coordinates $\ldots,q^{i},\ldots,P_{i},\ldots$, the Hamilton equations on
$(M\times Q,\Omega_{\mathcal{F}},H_{\mathcal{F}})$ look particularly simple
\begin{align*}
\dot{q}^{i} &  =0\\
\dot{P}_{i} &  =-\frac{\partial H_{\mathcal{F}}}{\partial q^{i}}%
=\mathrm{const.}%
\end{align*}
In this sense $W$ \emph{generated a canonical transformation that simplifies
the original problem}. In this respect, see \cite{mmm09}, where a quantum
version of the last result is also proposed. It is an interesting issue
developing Hamilton-Jacobi techniques for the computation of quantum
propagators (see, e.g., \cite{w41}).

\subsubsection{Hamiltonian dynamics of boundary data}

In the case when the Hamiltonian system $(T^{\ast}M,\Omega,H)$ comes from a
regular Lagrangian system, with Lagrangian $L\in C^{\infty}(TM)$, then one can
choose $Q=M$ and there is a canonical choice for $W$, namely
\[
W(x,q)=\int_{t_{0}}^{t_{1}}L(\gamma(t),\dot{\gamma}(t))dt
\]
where $\gamma$ is the solution of the Euler-Lagrange equations such that
$\gamma(t_{0})=x$, and $\gamma(t_{1})=q$. In this case, the diffeomorphism
$\Phi^{-1}:T^{\ast}M\simeq M\times Q$ \textquotedblleft\emph{transforms the
symplectic manifold of initial data of Hamilton equations into a symplectic
manifold of boundary data of the Euler-Lagrange equations}\textquotedblright.
In particular, there is a Hamiltonian system on boundary data (see
\cite{mmm09}). Rovelli \cite{r03} showed that this considerations can be
generalized (in a covariant way) to any classical field theory. In particular
he was able to write down a Hamilton-Jacobi equation on the space of boundary
data of the field equations and show that the action functional provide a
canonical solution. He also used the Hamilton-Jacobi equation to perform a
\emph{transition to the quantum regime}. In the case of Einstein gravity, he
obtained the Wheeler-De Witt equation (see also \cite{b66}). Unfortunately,
Rovelli's theory is rather far from being fully general and mathematically
rigorous. More recently, I and G.{} Moreno proved \cite{mv13} that whatever
the Lagrangian field theory one starts from (any number of dependent and
independent variables, derivatives, and gauge symmetries), \emph{there is a
canonical Hamiltonian system on the space of boundary data} (see also
\cite{m12}) which is, in a sense, equivalent to the Euler-Lagrange equations.
We achieved this result in full rigour within the jet space (and, in
particular, $\infty$-jet space) approach to PDEs. However, it is not clear
what is the precise relation to Hamilton-Jacobi theory. In particular, it is
not clear in what sense the action provide a complete integral of the
Hamilton-Jacobi problem, nor if one could actually quantize along this lines.
Notice that, in this formalism, characteristic Cauchy data, and singularities
of solutions should play a distinguished role. Clarifying these issues is, in
my opinion, an interesting open problem.

\section*{Conclusions}

Consider a determined system of quasi-linear PDEs governing the dynamics of a
field in the space-time. The boundary of a disturbance in the field, i.e., a
wave-front, is a characteristic surface in the space-time. In their turn,
wave-fronts propagate along bicharacteristics and bicharacteristics are often
trajectories of a Hamiltonian system. I just described a mathematically rather
precise way to pass from waves to rays, or from fields to particles. The
possibility of making this passage may be understood as a manifestation of the
quantum-mechanical wave-particle duality. Accordingly, the transition
\begin{equation}
\text{field equations }\Longrightarrow\text{ characteristic surfaces
}\Longrightarrow\text{ bicharacteristics}\label{10}%
\end{equation}
may be understood as analogous to the transition
\begin{equation}
\text{quantum mechanics }\Longrightarrow\text{ short wave-lenght limit
}\Longrightarrow\text{ classical mechanics.}\label{11}%
\end{equation}
Notice that in both transitions one progressively lose information. Therefore,
it should be expected that, performing the inverse transitions (in particular,
\emph{quantizing}) requires additional information. Vinogradov conjectured
(see the appedix of \cite{b...99}) that part of this information is actually
contained in the singularity equations. The idea that the geometric theory of
PDEs can account for quantization is intriguing and worth to be explored.

\subsection*{Acknowledgments}

I thank the scientific committee of the XXII International Fall Workshop on
Geometry and Physics for the opportunity to give the mini-course on which this
review is based. I also thank Beppe Marmo for suggesting me the topic of the
mini-course, for carefully reading a preliminary version of this paper, and
for all its valuable comments.

\end{document}